\documentclass[a4paper]{amsart}

\usepackage[latin1]{inputenc}
\usepackage{epsfig}
\usepackage{amsmath, amsthm, amssymb, graphicx}
\usepackage[all]{xy}

\newtheorem{theorem}{Theorem}[section]
\newtheorem{lemma}[theorem]{Lemma}

\newtheorem{prop}[theorem]{Proposition}

\newtheorem{cor}[theorem]{Corollary}
\newtheorem{thm}{Theorem}
\theoremstyle{definition}
\newtheorem{definition}[theorem]{Definition}

\def\co{\colon\thinspace}
\DeclareMathOperator{\fin}{fin}
\DeclareMathOperator{\kdef}{def}
\DeclareMathOperator{\id}{id}
\DeclareMathOperator{\Int}{Int}

\title{\bf Angle-deformations in Coxeter groups}

\author{Timothée Marquis}
\address{D\'epartement de Math\'ematiques \\
Université Libre de Bruxelles\\\newline
Boulevard du Triomphe CP 216\\
1050 Brussels\\
Belgium}
\email{tmarquis@ulb.ac.be}

\author{Bernhard Mühlherr}
\address{Mathematisches Institut \\
Universität Giessen\\\newline
Arndtstrasse 2\\
D-35392 Giessen\\
Germany}
\email{bmuhlher@ulb.ac.be}
\email{Bernhard.Muehlherr@math.uni-giessen.de}

\begin{document}

\begin{abstract}
The isomorphism problem for Coxeter groups has been 
reduced to its 'reflection preserving version' by 
B. Howlett and the second author. Thus, in order 
to solve it, it suffices to determine for a given Coxeter system 
$(W,R)$ all Coxeter generating sets $S$ of $W$ which 
are contained in $R^W$, the set of reflections of $(W,R)$. In this 
paper, we provide a further reduction: 
it suffices to determine all Coxeter generating sets $S \subseteq R^W$ 
which are sharp-angled with respect to $R$. 
\end{abstract}

\maketitle

\section{Introduction}

Let $W$ be a group and let $R \subseteq W$.
We call $R$ a \emph{Coxeter generating set}
of $W$ if $(W,R)$ is a Coxeter system.
All Coxeter systems $(W,R)$ considered in this paper
are assumed to have finite rank, i.e. $R$ is a finite
set.

Let $(W,R)$ be a Coxeter system and let
$S \subseteq R^W$ be
a Coxeter generating set
of $W$. A subset $J$ of $S$ is called \emph{spherical}
if it generates a finite subgroup; if it is of cardinality
2, it is called an 
\emph{edge} of $S$. Let  
$\{ s,t \}\subseteq S$ be an edge of $S$.  
By basic results on 
Coxeter groups, one knows that there exist
$r,r' \in R$ and $w \in W$ such that
$\langle s,t \rangle^w = \langle r,r' \rangle$.
If there exist $r,r' \in R$ and $w \in W$
such that $\{ s,t \}^w = \{ r,r' \}$, then
we call the edge $\{ s,t \}$ \emph{sharp-angled}
with respect to $R$. We call $S$ \emph{sharp-angled}
with respect to $R$ if all edges of
$S$ are sharp-angled with respect to $R$.  
The trivial example
of the dihedral groups shows that there
are examples of Coxeter systems $(W,R)$ admitting
Coxeter generating sets $S \subseteq R^W$
which are not sharp-angled with respect to $R$.

In Mühlherr \cite{toro}, it was conjectured that for any
Coxeter generating set $S \subseteq R^W$, there exists
an automorphism $\alpha$ of $W$ such that 
$\alpha(S) \subseteq R^W$ and such that $\alpha(S)$
is sharp-angled (Conjecture 1 in loc. cit.).
This conjecture may be seen as a reduction step
in order to state the main conjecture about
the solution of the isomorphism problem for
Coxeter groups, which is Conjecture 2 in \cite{toro}
(see Remark 1 below).
  
It was mentioned without proof in \cite{toro} 
that Conjecture 1 is true if there is no subdiagram
of type $H_3$. It turned out that this conjecture
was too optimistic if there are $H_3$-subdiagrams. 
Counter-examples 
have been found independently
by Ratcliffe and Tschantz and by Grassi (see \cite{RT} and \cite{Gr}).
This motivates the question whether it is still true
that one can reduce the solution for the isomorphism problem
to the main conjecture.
The goal of this paper is to show that this is indeed 
the case.

Our first result is the following.

\begin{thm}\label{thm1}
Let $(W,R)$ be a Coxeter system. 
Let $S \subseteq R^W$ be a Coxeter generating set of $W$
having no subsystem of type $H_3$.
Then there exists an automorphism $\alpha$ of $W$
such that $\alpha(S)$ is sharp-angled with respect to $R$.
\end{thm}

\noindent
As already mentioned before, Theorem \ref{thm1} has been announced
in \cite{toro} and it is a special
case of Theorem \ref{thm2} below. Its proof is given
in Section \ref{sec6}. We prefer to present it
separately since it is rather easy and provides
at the same time a good overview on the kind of
arguments that will yield Theorem \ref{thm2}.

The situation becomes considerably more complicated if $H_3$-subdiagrams are
allowed. First of all, the counter-examples to Conjecture 1
show that one cannot expect to produce sharp-angled Coxeter generating
sets from $S$ by automorphisms. So, we have to produce the 
desired Coxeter generating set starting
from $S$ by a sequence of operations which we call \emph{angle-deformations}.

In order to define \emph{angle-deformations}, we analyse the situation
where we are given a Coxeter system $(W,R)$ and a Coxeter generating
set $S  \subseteq R^W$ such that there is an edge  $J$ of $S$
which is not sharp-angled with respect to $R$. It turns out
that the Coxeter diagram of the system $(W,S)$ has to satisfy
several conditions with respect to
 the subset $J$. These conditions will be deduced
in Section \ref{figDE}. An edge satisfying these conditions
will be called a \emph{$\Delta$-edge}. 

Let $(W,S)$ be a Coxeter system and $J=\{ r,s \}$ be
a $\Delta$-edge of $S$.
Then we construct a mapping
$\delta\co S \to W$ such that 
$\delta(s) = s, \delta(r) \in \langle s,r \rangle$
and such that $S' := \{ \delta(x) \mid x \in S \}$
is a Coxeter generating set with the property that all spherical 2-subsets
$\{ x',y'  \} \neq \{ \delta(r),\delta(s) \}$ are sharp-angled
with respect to $S$. We call these mappings \emph{$J$-deformations}.
In the case where there are no $H_3$-subdiagrams, it is
easy to give the definition of these $J$-deformations.
If there are $H_3$-subdiagrams, the definition is given
recursively. We first define $J$-deformations for
a class of diagrams which we call {\it tame}.
The general case will then be treated by induction
on the number of `wild' vertices.

The construction of $J$-deformations will enable
us to prove our main result, which is the following.

\begin{thm}\label{thm2}
Let $(W,R)$ be a Coxeter system and
let $S \subseteq R^W$ be a Coxeter generating set of $W$.
Then, there exists a sequence $S = S_0,\ldots,S_k=S'$
of Coxeter generating sets $S_i$ such that
$S_i$ is a $J_i$-deformation of $S_{i-1}$
for some $\Delta$-edge $J_i$ of $S_{i-1}$ for each $1 \leq i \leq k$,
and such that $S'$ is sharp-angled with respect to $R$. 
\end{thm}

\noindent
We remark that the proof of Theorem \ref{thm2} is constructive.
Hence it provides 
a concrete algorithm to obtain 
the set $S'$ starting from $S$. Combining the theorem above
with the fact that the isomorphism problem for Coxeter groups
is reduced to its `reflection-preserving version' (as described
in \cite{toro}), we obtain the following.

\begin{cor}  
The isomorphism problem for Coxeter groups is solved as soon
as the following problem is solved.
\end{cor}

\noindent
{\bf Problem:} Let $(W,R)$ be a Coxeter system. Find 
all Coxeter generating sets $S \subset R^W$ such that
$S$ is sharp-angled with respect to $R$.

\subsection*{Remarks}

\noindent
{\bf 1.} There is a conjecture about the solution
of the above problem. This is Conjecture 2 in 
\cite{toro} and it is a refinement of Conjecture~8.\,1
in Brady--McCammond--Mühlherr--Neumann \cite{BMMN}. It says that if $R$ and $S$ are
as in the problem above, one can transform $S$ into
$R$ by a sequence of twists which had been introduced
in \cite{BMMN}. 
The conjecture has been proved for various
classes of Coxeter systems; the reader may refer to 
\cite{toro} for a survey on its status in 2005.
Recently, it was shown by Ratcliffe and Tschantz
in \cite{RT} that the conjecture holds for
chordal Coxeter systems as well.

\noindent
{\bf 2.}
In \cite{RT}, our main result has been obtained for
chordal Coxeter systems. Their methods are quite
different from ours. Their arguments rely heavily
on a very strong property of chordal Coxeter groups
which is not available in the general case.

\smallskip
\noindent
The paper is organized as follows. In Section \ref{sec2},
we fix notation, recall some basic facts on Coxeter
groups and provide some preliminary results. 
In Section  \ref{sec3}, we introduce angle-deformations
and make some observations about them. In Section \ref{sec4},
we prepare the proof of Theorem \ref{thm1}. In this section
we introduce $\Theta$-edges, which are special cases
of $\Delta$-edges.
Section \ref{sec5} is devoted to introduce and
investigate the notion of a sharp-angled set
of reflections in a Coxeter group. This will 
enable us to give the proof of  Theorem \ref{thm1} in Section 
\ref{sec6}. In Section \ref{sec7}, we collect several 
informations about angle-deformations of Coxeter systems
with subdiagrams of type $H_3$ and $H_4$. In Section \ref{figDE},
we define $\Delta$-edges. Later on, these turn out to be
precisely the edges of a Coxeter system for which there
are non-trivial angle-deformations. This fact
is a consequence of Proposition \ref{prfthm1prop} and
Theorem \ref{mainthm}, and it is
indeed the key-ingredient of the proof of our main result.
Section \ref{sec9} can be seen as a preparatory section
for the proof of Theorem \ref{mainthm} which will be 
completed in Section \ref{tausection}. In Section \ref{sec11},
we finally give the proof of our main result.

\section{Preliminaries} 
\label{sec2}

\subsection{Graphs}
For a set $X$, denote by $P_2(X)$ the set of all subsets of $X$ having cardinality 2. A \emph{graph} is a pair $(V,E)$ consisting of a set $V$ and a set $E\subseteq P_2(V)$. The elements of $V$ and $E$ are called \emph{vertices} and \emph{edges} respectively.

Let $\Gamma=(V,E)$ be a graph. Let $v$,$w$ be two vertices of $\Gamma$. They are called \emph{adjacent} if $\{v,w\}\in E$. In this paper, a \emph{path} from $v$ to $w$ is a sequence $v=v_0, v_1, \dots, v_k=w$, where $v_{i-1}$ is adjacent to $v_i$ for all $1\leq i\leq k$ and where $v_1, \dots, v_k$ are pairwise distinct; the number $k$ is the \emph{length} of the path. The path is \emph{minimal} if it is of minimal length. The \emph{distance} between $v$ and $w$ (denoted by $\delta(v,w)$) is the length of a minimal path joining them; if there is no path joining $v$ and $w$, we put $\delta(v,w)=\infty$.

A path $v=v_0, v_1, \dots, v_k=w$ is said to be \emph{chordfree} if $E\cap P_2(\{v_0, \dots, v_k\}) = \{\{v_0,v_1\}, \{v_1,v_2\}, \dots, \{v_{k-1},v_k\}\}$. A path $v=v_0, v_1, \dots, v_k=w$ is called a \emph{circuit} if $v=w$ and $k\geq 2$.

The relation $R\subseteq V\times V$ defined by $R=\{(v,w)|\delta(v,w)\neq \infty\}$ is an equivalence relation whose equivalence classes are called the \emph{connected components} of $\Gamma$. A graph is said to be \emph{connected} if it has only one connected component.

\subsection{Coxeter systems}
Let $(W,S)$ be a pair consisting of a group $W$ and a set $S\subseteq W$ of involutions. For $r,s\in S$, denote by $m_{rs}\in \mathbb{N}\cup \{\infty\}$ the order of the product $rs$ in $W$. Note that we will also use the notation $o(rs)$ instead of $m_{rs}$. Define $E(S):= \{\{r,s\}\subseteq S \ | \ 1\neq m_{rs}\neq \infty\}$ to be the set of \emph{edges} of $S$. Then $\Gamma(S)$ is the graph $(S,E(S))$ whose edges are labelled by the corresponding $m_{rs}$. Throughout this text, any graph notion (such as paths and circuits) associated to the pair $(W,S)$ must be understood as being in $\Gamma(S)$. In particular, when we speak about the "diagram of $(W,S)$", we refer to $\Gamma(S)$.

The \emph{Coxeter diagram} associated to $(W,S)$ is the graph $(S,E'(S))$ where $E'(S):= \{\{r,s\}\subseteq S | m_{rs}\geq 3\}$ and where the edges are labelled by the corresponding $m_{rs}$. A subset $K$ of $S$ is said to be \emph{irreducible} if the underlying Coxeter subdiagram $(K,E'(K))$ is connected. We call $K$ \emph{spherical} if it generates a finite group. Finally, $K$ is \emph{2-spherical} if $m_{rs}< \infty$ for all $r,s\in K$. If $S$ is irreducible, spherical or 2-spherical, we say that $(W,S)$ is \emph{irreducible}, \emph{spherical} or \emph{2-spherical}, respectively. Note that sometimes, we use the same notions for $\Gamma(S)$ instead of $(W,S)$.

We say that $(W,S)$ is a \emph{Coxeter system} if $S$ generates $W$ and if the relations $((rs)^{m_{rs}})_{r,s\in S}$ form a presentation of $W$. We call $R\subseteq W$ a \emph{Coxeter generating set} if $(W,R)$ is a Coxeter system. 

Let $(W,R)$ be a Coxeter system. An element of $W$ is called a \emph{reflection} if it is conjugate in $W$ to an element of $R$; the set of all reflections is denoted by $R^W$.

\subsection{Conventions about figures}
Here are some conventions about the figures appearing in the paper, which the reader may refer back to when needed.
\\
Throughout this text, all figures represent diagrams of the form $\Gamma(K)$ for some Coxeter system $(W,S)$ and $K\subseteq S$. The edges in plain have a finite label, while the edges in strips have an infinite label. An absence of edge does not imply anything. If there is a single edge with more than one label (say $m>1$), then the figure must be understood as $m$ different figures, one for each of these labels. If there are two or more edges with more than one label, then these edges will have the same number $m>1$ of labels. In that case, the figure must be understood as $m$ different figures, the $i$-th figure being obtained by taking the $i$-th label from each of these edges.
\\
A dotted line between two vertices means that there is a path (in plain) joining these two vertices, but the other vertices in the path were omitted. (It will be always clear from the context what the omitted vertices are). For example, in Section \ref{figDE}, figures \ref{figDE3} and \ref{figDE4} contain a path $\{S(1),S(2),\dots,S(n-1),S(n)\}$. We denote by $X$ this set and we assume $n\geq 2$. Let $X_1:=X\setminus \{S(1)\}$ and $X_n:=X\setminus \{S(n)\}$. We assume $X$ has the following property:
$$o(S(i)S(j))= \infty \ \textrm{for all $i,j$ such that} \ 1\leq i < j \leq n \ \textrm{and} \ |i-j| \geq 2.$$
Finally, for a vertex $y\notin X_1$, we mean by $\overline{X_1y}=\infty$ that $m_{xy}=\infty$ for all $x\in X_1$.

\subsection{Coxeter generating sets and automorphisms}

\begin{lemma}\label{amallemma}
Let $(W,S)$ be a Coxeter system and let $S_1, S_2$ be subsets of $S$ such that each edge of $S$ is contained in $S_1$ or $S_2$. Put $S_0 := S_1 \cap S_2$. Let $\delta\co S \to W$ be a mapping such that $\delta(S_i)$ is a Coxeter generating set of $\langle S_i \rangle$ for $i = 0,1,2$. Then $\delta(S)$ is a Coxeter generating set of $W$. Moreover, if the restrictions of $\delta$ to $S_1$ and $S_2$ extend to automorphisms of $\langle S_1 \rangle$ and $\langle S_2 \rangle$ respectively, then $\delta$ extends to an automorphism of $W$.
\end{lemma}

\begin{proof}
 This follows immediately from the fact that
$W = \langle S_1 \rangle *_{\langle S_0 \rangle} \langle S_2 \rangle$. 
\end{proof}

\noindent
The following lemma follows easily by the pigeon-hole principle.

\begin{lemma} \label{piwhlemma}
Let $G$ be a finite group, let $\alpha$ be an automorphim of $G$
and let $g \in G$. Then $\alpha^m(g)\alpha^{m-1}(g)\ldots\alpha^2(g)\alpha(g)g 
= 1_G$ for some $m \geq 0$.
\end{lemma}

\noindent
Using the previous lemma, one immediately obtains the 
following proposition.

\begin{prop} \label{autoprop}
Let $(W,S)$ be a Coxeter system and let 
$\alpha\co W \to W$ be an epimorphism. Suppose that
there is a subset ${\bf K}$ of $2^S$ such that the following holds:
\begin{itemize}
\item[(1)] All elements of ${\bf K}$ are spherical.
\item[(2)] For all $K \in {\bf K}$, the mapping
$\alpha \mid_{\langle K \rangle}$ is an automorphism of 
$\langle K \rangle$.
\item[(3)] For all $s \in S$, there exists 
$w_s \in \bigcup_{K \in {\bf K}} \langle K \rangle$ such that
$\alpha(s) = w_s s w_s^{-1}$.
\end{itemize}
Then $\alpha$ is an automorphism of $W$ which is of finite order.
\end{prop}

\subsection{The geometric representation of a Coxeter system}

In this subsection, we collect several basic results about the geometric representation of a Coxeter system. The standard references are Bourbaki \cite{NB} and Humphreys \cite{JH}.
\\                   
Throughout this paper, $\Omega$ and $\Omega'$ are the following
subsets of $\mathbb{R}$: 
$$\Omega := \{ \cos(\pi/m) \mid m \in \mathbb{N} \} \cup [1,\infty)$$
and $\Omega' := \Omega \setminus \{-1 \}$.
Moreover, we define a mapping 
$C\co \mathbb{N} \cup \{ \infty \} \to -\Omega$
by setting $C(m) := -\cos(\pi/m)$ if $m \in {\bf N}$ and
$C(\infty) := -1$.

Let $V$ be a real vector space endowed with a symmetric bilinear form
$b\co V \times V \to \mathbb{R}$. The set of vectors $v \in V$
with $b(v,v)=1$ is denoted by $U(V,b)$ and for each
such vector, the corresponding orthogonal reflection
with respect to $b$ is denoted by $\rho_v$;
hence $\rho_v(x) = x -2b(x,v)v$ for each $x \in V$.

Let $(W,R)$ be a Coxeter system.  Let
$V := \mathbb{R}^R$ and  $(e_r)_{r \in R}$ be the canonical basis of $V$.
Furthermore, let $b\co V \times V \to \mathbb{R}$
be the symmetric bilinear form defined by
$b(e_r,e_s) := C(o(rs))$.
                                                                                
\begin{theorem} The mapping $r \mapsto \rho_{e_r}$
from $R$ into $O(V,b)$ extends to a monomorphism
from $W$ into $O(V,b)$.
\end{theorem}

\noindent
Thus, by the above construction, we obtain
a canonical faithful linear representation of the Coxeter
group $W$ which is called the \emph{geometric representation} of $(W,R)$.
We now identify $W$ with its image in $O(V,b)$ and we put
$\Phi(W,R) := \{ w(e_r) \mid w \in W, r \in R \}$. We have the following:

                                                                                \begin{lemma} \label{refllemma} For all $r \in R$ and $w \in W$, 
we have $\rho_{w(e_r)} = wrw^{-1}$; 
in particular, $R^W = \{ \rho_{\alpha} \mid \alpha \in \Phi(W,R) \}$.  
Moreover, if $\alpha,\beta \in \Phi(W,R)$ are such that 
$\rho_{\alpha} = \rho_{\beta}$, then $\beta = \alpha$ or $\beta = -\alpha$.
\end{lemma}

\noindent
The set $\Phi:=\Phi(W,R)$ is called the \emph{root system} of 
$(W,R)$.
We put $$V^+ := \{ \Sigma_{r \in R} \mu_r e_r \mid \mu_r \geq 0 \mbox{ for all }r \in R \}$$
and $V^- := -V^+$; furthermore, we put
$\Phi^+ := V^+ \cap \Phi$ and $\Phi^- := V^- \cap \Phi$.

\begin{lemma}
$\Phi = \Phi^+ \cup \Phi^-$.
\end{lemma}

\noindent
The elements of $\Phi^+$ are called the \emph{positive roots} of $(W,R)$.
A subset $\Pi$ of $\Phi$ is called a \emph{root-subbase} of $\Phi$
if $\Pi \subseteq \Phi^+$ and if
$b(\alpha,\beta) \in -\Omega'$ for all $\alpha \neq \beta \in \Pi$.

\noindent                        
The following theorem is a consequence of the main result in Deodhar \cite{VD2} and Dyer \cite{MD}.                                    
\begin{theorem} \label{Dyerthm}
Let $\Pi$ be a root subbase of $\Phi$
and put $S := \{ \rho_{\alpha} \mid \alpha \in \Pi \}$.
Then $(\langle S \rangle,S)$ is a Coxeter system.                                      
Conversely, let $W'$ be a subgroup of $W$ which is generated
by a set of reflections. Then there exists a root-subbase $\Pi'$
of $\Phi$ such that $W' = \langle \rho_{\alpha} \mid \alpha \in \Pi' \rangle$.
\end{theorem}

\subsection{Flexibility}

Let $(W,S)$ be a Coxeter system and $J\subseteq S$. We define the following notions and notations:
\begin{itemize}
\item
$J^{\bot} := \{s\in S \ | \ \forall \ j\in J \ : \ m_{sj}=2\}$.
\item
$J^{\fin} := \{s\in S\setminus J \ | \ m_{sj}<\infty \ \forall\ j\in J \}$.
\item
$J^{\infty} := \{s\in S\setminus J \ | \ \exists \ j\in J \ : \ m_{sj}=\infty\}$.
\item
$\mathcal{G}_J := (J^{\infty}, \ \{\{a,b\}\subseteq J^{\infty} \ | \ m_{ab}<\infty\})$.
\item
A \emph{$J$-component} is a connected component of $\mathcal{G}_J$.
\item
Let $L$ be a $J$-component. We shall say that $j\in J$ is \emph{$L$-free} if $m_{jl}=\infty$ for all $l\in L$.
\item
An element $j$ of $J$ is \emph{$J^{\infty}$-free} if it is $L$-free for every $J$-component $L$.
\item
The $J$-component $L$ is said to be \emph{flexible} if there exists $j\in J$ such that $j$ is $L$-free.
\item
Finally, we will say that $J$ is \emph{flexible} if all $J$-components are.
\end{itemize}
\noindent
Here is a first observation.

\begin{lemma} \label{flexedgelemma}
Let $(W,S)$ be a Coxeter system and let $J= \{ r,s \}$ be an edge of $S$.
Then $J$ is flexible if and only if there is no chordfree circuit in $\Gamma(S)$ of 
length at least 4 containing $J$.
\end{lemma}

\begin{figure}
\begin{center}
\includegraphics[trim = 0mm 100mm 0mm 0mm, clip, height=4cm]{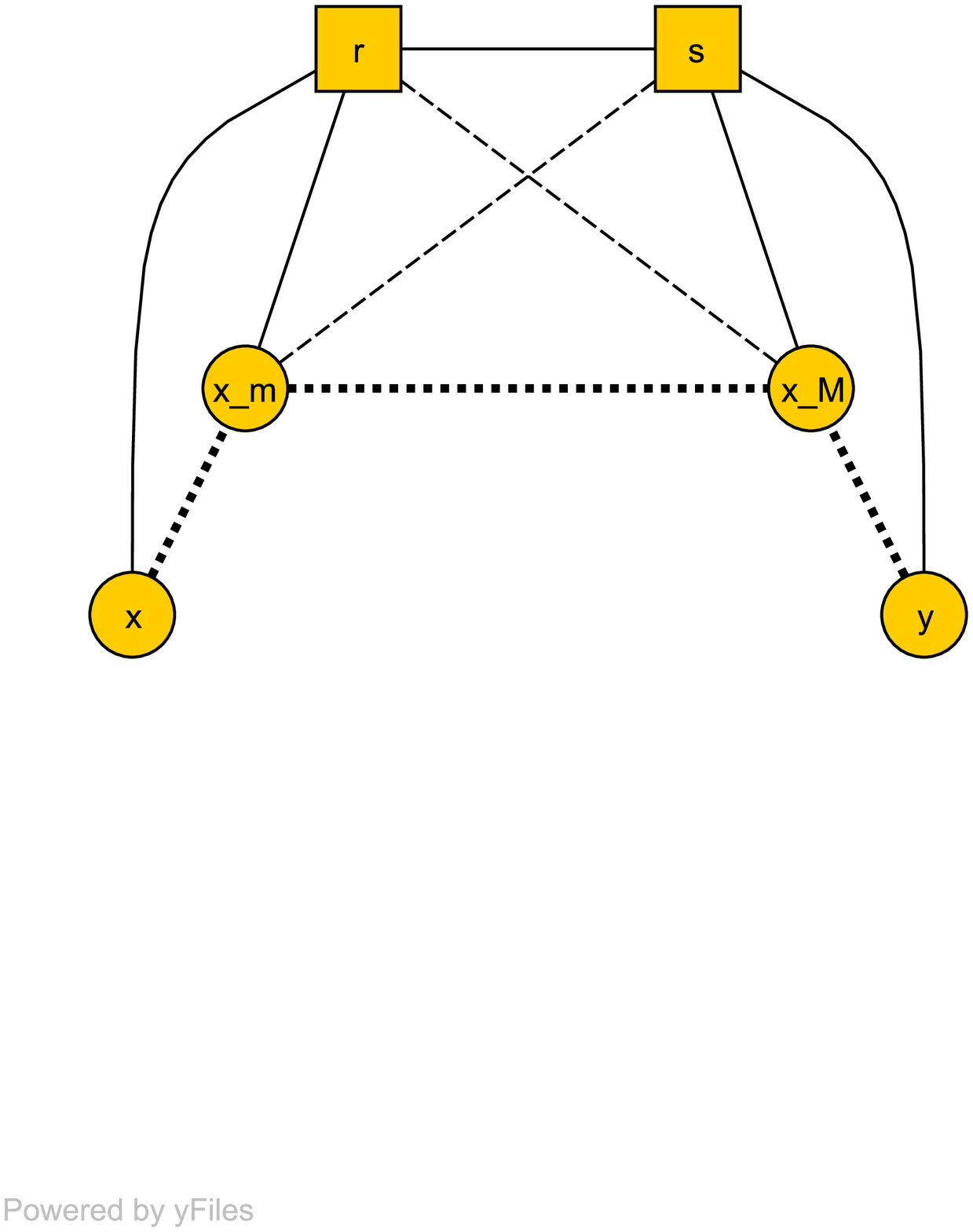}
\caption{Lemma \ref{flexedgelemma}.}
\label{graph1_Jflexible}
\end{center}
\end{figure}

\begin{proof}
Suppose first $J$ is not flexible. Then there exists a $J$-component $L$ and $x,y\in L$ such that $m_{xr}<\infty$ and $m_{ys}<\infty$.
Let $x=x_0, x_1, \dots, x_k = y$ be a minimal path in $L$ joining $x$ to $y$. Define
$$M := min\left\{i \ | \ 0 < i \leq k; \ m_{x_is}<\infty \right\}$$
and
$$m := max\left\{i \ | \ 0\leq i < M; \ m_{x_ir}<\infty \right\}.$$
Then the subpath $x_m, x_{m+1}, \dots , x_M$ from $x_m$ to $x_M$ is still minimal, hence chordfree, and possesses the following properties:
\begin{itemize}
\item[(1)]
$m_{x_is}=\infty$ for all $i$ such that $m\leq i<M$ (by definition of M);
\item[(2)]
$m_{x_ir}=\infty$ for all $i$ such that $m<i\leq M$ (by definition of m).
\end{itemize}
Moreover, $m_{x_Ms}<\infty$ and $m_{x_mr}<\infty$. We then obtain a chordfree circuit $r, x_m,$ $x_{m+1}, \dots, x_M, s, r$, as required. The situation is illustrated on figure \ref{graph1_Jflexible}. 
\\
The converse is obvious. 
\end{proof}

\section{Angle-deformations}
\label{sec3}

\begin{definition} Let $(W,S)$ be a Coxeter system,
let $J = \{ r,s \}$ be an edge of $S$ and let
$\omega \in \langle J \rangle$ be such that 
$\omega r \omega^{-1}$ and $s$ generate $\langle J \rangle$.
An \emph{$(r,s,\omega)$-deformation} of $S$ is a mapping
$\delta\co S \to W$ satisfying the following
properties:

\begin{itemize}
\item[AD1)] $\delta(x) \in S^W$ for all $x \in S$;
\item[AD2)] $\delta(r) = \omega r \omega^{-1}$ and 
$\delta(s) = s$;
\item[AD3)] $\delta(S)$ is a Coxeter generating set of $W$;
\item[AD4)] there exists a bijection $\Delta$ from the set of 
edges of $S$ onto the set of edges of $\delta(S)$ such that
$\Delta(J) = \{  \omega r \omega^{-1},s  \}$ and such that
for each edge $K \neq J$ of $S$, there exists $w_K \in W$
with $\Delta(K) = K^{w_K}$.
\end{itemize}

\end{definition}

\begin{definition} Let $(W,S)$ be a Coxeter system and
let $J = \{ r,s \}$ be an edge of $S$.
A \emph{$J$-deformation} of $S$ is an $(r,s,\omega)$-deformation
of $S$ for some $\omega \in \langle J \rangle$.
An \emph{angle-deformation} of $S$ is a $J$-deformation
for some edge $J$ of $S$.
\end{definition}

\noindent
The following Proposition is a consequence of Lemma \ref{amallemma}.
\begin{prop} \label{amalprop}
Let $(W,S)$ be a Coxeter system and let $S_1,S_2$ be subsets
of $S$ such that each edge of $S$ is contained in $S_1$ or $S_2$
and put $S_0 := S_1 \cap S_2$. Let $J$ be an edge contained
in $S_0$ and assume that $\delta_i\co S_i \to \langle S_i \rangle$
are $J$-deformations of $S_i$ for $i=0,1,2$ and that
$\delta_0 = \delta_i \mid_{S_0}$ for $i=1,2$.
Define $\delta\co S_1 \cup S_2 \to \langle S_1 \cup S_2 \rangle$
by setting $\delta \mid_{S_i} := \delta_i$ for $i=1,2$.
Then $\delta$ is a $J$-deformation of $S_1 \cup S_2$.
\end{prop}

\begin{prop} \label{defcritprop}
Let $(W,S)$ be a Coxeter system, $J:= \{ r,s \}$
be an edge of $S$ and let $\omega \in \langle J \rangle$ be such that
$\omega r \omega^{-1}$ and $s$ generate $\langle J \rangle$.
Let ${\bf K}$ be a set of spherical subsets
of $S$ such that each element of ${\bf K}$ contains $J$ and let
$\delta\co S \to W$ be a mapping with the following
properties:
\begin{itemize}
\item[a)] $\delta(r) = \omega r \omega^{-1}$ and $\delta(s)=s$;
\item[b)] $\langle \delta(S) \rangle = W$;
\item[c)] for all $x \in S$, there exists an element $w_x$
in $\bigcup_{K \in {\bf K}} \langle K \rangle$ such that
$\delta(x) = w_x x w_x^{-1}$;
\item[d)] for each edge $E$ of $S$ different from $J$, there exists
an element $w_E \in W$ such that $\delta(E) = E^{w_E}$.
\end{itemize}
Then $\delta$ is an $(r,s,\omega)$-deformation of $S$ 
which extends to an automorphism of $W$.
\end{prop}

\begin{proof} By the universal property of $(W,S)$ and
Property d), $\delta$ extends to an endomorphism $\alpha$ of $W$
which is in fact an epimorphism because of Property b).
By Proposition \ref{autoprop}, it follows now from
Property c) that $\alpha$ is an automorphism. Hence
$\delta(S)$ is a Coxeter generating set of $W$ and
the mapping $E \mapsto \delta(E)$ is a bijection as required
in Condition AD4). As AD1) is a consequence of Property c),
and as AD2) is precisely Property a), the proposition is proved.
\end{proof}

\section{Angle-deformations involving dihedral groups}
\label{sec4}

Throughout this section, $(W,S)$ is a Coxeter system
and $J = \{ r,s \}$ is an edge of $S$ such that $o(rs) \geq 3$.

\subsection{Condition (TWa)}

\begin{definition} Let $a \in J$. We say that $J$ is an \emph{$a$-special}
subset of $S$ if the following condition (TWa) holds.

\begin{itemize}
\item[(TWa)] For all $x \in S \setminus J$ we have $o(xa) \in
\{ 2,\infty \}$, and if $o(xa) =2$ then $x \in J^{\perp}$.
\end{itemize}
\end{definition}

\noindent
The following observation is immediate.

\begin{lemma} \label{lemma8.1} 
Let $a \in J$ be such that $J$ is $a$-special.
Then the following holds.
\begin{itemize}
\item[a)] $\{ J, J^{\infty},  J^{\perp} \}$ is a partition of $S$;
\item[b)] $a$ is $J^{\infty}$-free; in particular, $J$ is fexible.
\end{itemize}
\end{lemma}

\noindent
For the remainder of this subsection, we assume that 
$a \in J$ is such that $J$ is $a$-special, and
$\omega \in \langle J \rangle$ is such that 
$\omega r \omega^{-1}$ and $s$ generate $\langle J \rangle$.
We put $\pi := 1_W$ if $a = r$ and $\pi := \omega$ if $a=s$. 
Moreover, we let $\delta\co S \to W$ be the
mapping defined by $\delta(r) = \omega r \omega^{-1}$,
$\delta(y) = y$ for $y \in \{ s \} \cup J^{\perp}$
and $\delta(x) := \pi x \pi^{-1}$ if $x \in J^{\infty}$.

\begin{lemma} \label{lemma8.2} 
Let $E = \{ x,y \}$ be an edge of $S$
different from $J$.
Then there exists $w_E$ such that
$\delta(E) = E^{w_E}$.
\end{lemma}

\begin{proof} 
Note first that each $y \in J^{\perp}$ commutes with $\omega$ and $\pi$.
Hence, if $E \subseteq \{ s \} \cup J^{\perp}$, then we may choose $w_E = 1_W$;
if $E \subseteq \{ r \} \cup J^{\perp}$, then we may choose $w_E = \omega$;
and if $E \subseteq J^{\infty} \cup J^{\perp}$, then we may choose $w_E = \pi$.

By the previous lemma, we 
are left with the case where $E \subseteq J \cup J^{\infty}$.
As $a$ is $J^{\infty}$-free and $E \neq J$, we are now left with
the case where $E = \{ b,x \}$ for some $x \in J^{\infty}$
and where $b$ is the element of $J$ distinct from $a$.
If $a = r$, we may choose $w_E =  1_W$ and if $a = s$, we may choose
$w_E = \omega$.
\end{proof}

\begin{prop} \label{rk2cpdefprop} 
The mapping $\delta$ is an $(r,s,\omega)$-deformation
of $S$ which extends to an automorphism of $W$.
\end{prop}

\begin{proof} Setting ${\bf K} := \{ J \}$, Properties a), b) and c)
required in Proposition \ref{defcritprop} are clear from the definition
of $\delta$ and Property d) is settled by the previous lemma. 
\end{proof}

\subsection{$\Theta$-edges}

\begin{definition} We say that $J$ is a \emph{$\Theta$-edge} of $S$
if $J$ is flexible and if there is no 2-spherical and irreducible
subset of $S$ containing $J$ properly.
\end{definition}

\noindent
{\bf Remark:} If $J$ is a $\Theta$-edge, then $\{ J,J^{\infty},J^{\perp} \}$
is a partition of $S$.

\noindent
For the remainder of this subsection, we suppose that $J$ is a $\Theta$-edge of $S$.
Moreover, we assume that $\omega \in \langle J \rangle$ is such that
$\omega r \omega^{-1}$ and $s$ generate $\langle J \rangle$.

\noindent
Let $L$ be a $J$-component.
We denote the set of $L$-free vertices in
$J$ by $\Pi(L)$. It is non-empty because $J$ is assumed to be flexible.
If  $r \in \Pi(L)$, we put $a_L := r$ and $\gamma_L := 1_W$; if this is not
the case, we set  $a_L := s$ and $\gamma_L := \omega$. 
We set $K_L := J \cup L \cup J^{\perp}$. 
We define $\delta_L\co K_L \to \langle K_L \rangle$ by
$\delta_L(r) := \omega r \omega^{-1}$, $\delta_L(y) :=y$ for all $y \in \{ s \} \cup J^{\perp}$
and $\delta_L(x) := \gamma_L x \gamma_L^{-1}$ for all $x\in L$.

\begin{prop} \label{rk2defprop}
Let $\delta\co S \to W$ be the unique mapping such that $\delta \mid_{K_L} = \delta_L$ for every $J$-component $L$.
Then $\delta$ is an $(r,s,\omega)$-deformation of $S$ wich extends to an automorphism
of $W$.
\end{prop}

\begin{proof} Let $L$ be a $J$-component. The edge $J$ is an $a_L$-special 
subset of $K_L$ and hence it follows by Proposition
\ref{rk2cpdefprop} that $\delta_L$ is an $(r,s,\omega)$-deformation of
$K_L$. An obvious induction on the number of $J$-components
using Proposition \ref{amalprop} and Lemma \ref{amallemma} yields the claim.
\end{proof}

\section{Sharp-angled sets of reflections}
\label{sec5}

Throughout this section, $(W,R)$ denotes a Coxeter system,
where $W$ is identified with its image in $O(V,b)$
by its geometric representation and $\Phi \subseteq U(V,b)$
is its root system.

\begin{lemma} \label{howlettlemma}
Let $\alpha,\beta \in \Phi$.
\begin{itemize}
\item[a)] If 
$|b(\alpha,\beta)|<1$, then $o(\rho_{\alpha}\rho_{\beta})$
is finite and $b(\alpha,\beta) = -\cos(p\pi/q)$
for some integers $p$ and $q$. 
\item[b)] If $\rho_{\alpha} \neq \rho_{\beta}$
and $|b(\alpha,\beta)|\geq 1$, then $\rho_{\alpha}\rho_{\beta}$
has infinite order.
\item[c)] If $\rho_{\alpha} \neq \rho_{\beta}$, then 
$o(\rho_{\alpha}\rho_{\beta})$ is finite if and only
if $|b(\alpha,\beta)|<1$.
\end{itemize}

\end{lemma}
                             
\begin{proof} Assertion a) is Proposition 1.4 in Brink and Howlett \cite{BH93},
whereas Assertion b) is an easy exercice in linear algebra.
Assertion c) is an immediate consequence of a) and b).
\end{proof}                                                   
                                                                                
\begin{definition} Let $s \neq t \in R^W$ be such
that $o(st)$ is finite.
Let $\alpha,\beta \in \Phi$ be such that
$s = \rho_{\alpha}$ and $t = \rho_{\beta}$.
Then we call the 2-set $\{ s,t \}$ \emph{sharp-angled}
if $|b(\alpha,\beta)| \in \Omega$.
\end{definition}

\noindent
{\bf Remark:} Note that this definition does not depend
on the choice of $\alpha$ and $\beta$ in view of
the last statement of Lemma \ref{refllemma}.

\noindent
The following  two lemmas are easy.

\begin{lemma} \label{2346lemma}
Let $s \neq r \in R^W$ be such
that $o(rs)$ is finite.
If $\{ r,s \}$ is not sharp-angled, then
$o(rs) \geq 5$.
\end{lemma}

\begin{lemma} \label{rk2deflem}
Let $s \neq r \in R^W$ be such
that $o(rs)$ is finite and suppose $\{r,s\}$ is not sharp-angled.
Then there exists an element $w \in \langle s,r \rangle$
such that the set $\{ s,wrw^{-1} \}$ is sharp-angled. Moreover, if $o(rs)=5$, we may choose $w$ to be $srs$.
\end{lemma}

\begin{definition} A set $S \subset  R^W$ is called
\emph{sharp-angled} if each edge of $S$ is sharp-angled.
\end{definition}

\noindent
The following lemma follows from the fact that
$W$ is a subgroup of $O(V,b)$ and from the first statement of
Lemma \ref{refllemma}.
                                                                                
\begin{lemma} \label{sharpconjlem}
Let $S$ be a set of reflections and let
$w \in W$. Then $S^w$ is sharp-angled if and only if
$S$ is sharp-angled.
\end{lemma}

\noindent
The following fact follows from the definition of a root-subbase:
                                                                                
\begin{lemma} \label{Subbaselemma}
Let $\Pi$ be a root-subbase of $\Phi$
and $S:= \{ \rho_{\alpha} \mid \alpha \in \Pi \}$.
Then $S$ is sharp-angled.
\end{lemma}

\subsection{Fundamental sets of reflections}

\begin{definition} A subset $S$ of $R^W$ is called \emph{fundamental}
if $(\langle S \rangle,S)$ is a Coxeter system.
\end{definition}
                                                                        
\begin{theorem} \label{reflrigidthm}
Let $S \subset R^W$ be a fundamental set of
reflections and suppose that one of the following holds:
\begin{itemize}
\item[A)] The Coxeter system $(\langle S \rangle,S)$
is 2-spherical, irreducible and non-spherical.
\item[B)] $\Gamma(S)$ is a chordfree circuit
of length at least 4.
\end{itemize}
Then $S$
is sharp-angled.
\end{theorem}

\begin{proof} As $W':= \langle S \rangle$ is generated by
a set of reflections, we may apply the second part
of Theorem \ref{Dyerthm} to see that
there is a root-subbase $\Pi$ of $\Phi$
such that the set $S' := \{ \rho_{\alpha} \mid \alpha \in \Pi \}$
is a Coxeter generating set of $W'$.
It is known by Caprace--Mühlherr \cite{CM} and Charney--Davis \cite{CD}
that the Coxeter system $(\langle S \rangle,S)$
is strongly reflection rigid and hence $S$ and $S'$ are conjugate
in $W'$ and the claim follows from Lemmas \ref{sharpconjlem} and \ref{Subbaselemma}.
\end{proof}
                                                                                
\section{Proof of Theorem \ref{thm1}}
\label{sec6}

Throughout this section, $(W,R)$ is a Coxeter system and
$S \subseteq R^W$ is a fundamental set of reflections.
Moreover, we assume that $S$ contains no subset of
type $H_3$.

\begin{prop} \label{prfthm1prop}
Suppose that $J$ is an edge of $S$ which is not sharp-angled.
Then $J$ is a $\Theta$-edge of $S$.
\end{prop}

\begin{proof} Put $J = \{ r,s \}$.
By Lemma \ref{2346lemma}, we have $o(rs) \geq 5$. 
Let $t \in S$ be such that $o(rt)$ and $o(st)$ are finite. 
By Theorem \ref{reflrigidthm} and our hypothesis that there
are no subsets of type $H_3$, we have that $t\in J^{\perp}$.
Hence there is no irreducible 2-spherical subset of $S$
containing $J$ properly.
Furthermore, again by Theorem \ref{reflrigidthm},
 there is no chordfree circuit of length
at least 4 containing $J$. By Lemma \ref{flexedgelemma}, it follows
that $J$ is flexible. Hence $J$ is indeed a $\Theta$-edge of $S$.
\end{proof}

\begin{cor} \label{prfthm1cor}
Suppose that $J$ is an edge of $S$ which is not sharp-angled.
Then there exists a $J$-deformation $\delta$  of $S$
such that $\delta(J)$ is sharp-angled and such that $\delta$
is the restriction of an automorphism of $\langle S \rangle$.
\end{cor}

\begin{proof} Put $J = \{ r,s \}$. 
By Lemma \ref{rk2deflem}, we can find an element $\omega \in \langle J \rangle$
such that $\omega r \omega^{-1}$ and $s$ generate $\omega \in \langle J \rangle$
and such that $ \{ \omega r \omega^{-1}, s \}$ is sharp-angled.
By the previous proposition, we know that $J$ is a $\Theta$-edge of $S$
and hence, by Proposition \ref{rk2defprop}, we can find an $(r,s,\omega)$-deformation of $S$ which extends to an automorphism of $\langle S \rangle$ 
and we are done.
\end{proof}

\subsection*{Conclusion of the Proof of Theorem \ref{thm1}}

Let $S \subset R^W$ be a Coxeter generating set which is not sharp-angled.
Suppose $S$ contains $n \geq 1$ edges which are not sharp-angled and choose one of them. Call it $J$. By the previous corollary, there exists
a $J$-deformation $\delta$ of $S$ which extends to an automorphism of $W$
(because $\langle S \rangle = W$) and such that $\delta(J)$ is sharp-angled.
Let $J'$ be an edge of $S$ different from $J$. Then $\delta(J')$ is $W$-conjugate
to $J'$ by Property AD4) of $\delta$; in particular, $\delta(J')$ is sharp-angled
if and only if $J'$ is sharp-angled. Hence the number of edges in $\delta(S)$
which are not sharp-angled is $n-1$. Thus the statement follows by an
obvious induction on the number of edges of $S$
which are not sharp-angled. \qedsymbol

\section{Angle-deformations involving $H_k$}
\label{sec7}

\subsection{Coxeter systems of type $H_3$}

\begin{lemma} \label{H3basiclem}
Let $(W,S)$ be a Coxeter system of type $H_3$, where $S = \{ r,s,t \}$
and $o(rs)=5$, $o(st)=3$. Set $\omega := tsrtst$, $\pi := trs$
and define $\delta\co S \to W$ by
$\delta(r):= rsr,\delta(s):=s$ and $\delta(t):=\omega t \omega^{-1}$.
Then we have the following:

\begin{itemize}
\item[(1)] $\omega s \omega^{-1} = s$, $\omega t \omega^{-1} = \pi r \pi^{-1}$, $\pi t \pi^{-1}=rsr$.
\item[(2)] There is an automorphism $\alpha$ of $W$ which extends $\delta$.
\item[(3)] $\delta$ is an $(r,s,srs)$-deformation of $(W,S)$.
\end{itemize}

\end{lemma}

\begin{proof} Part (1) is a straightforward calculation. Moreover, it is clear 
that $\delta(S)$ is contained in $S^W$ and that it generates $W$.
It follows from (1) that $\{ \delta(s), \delta(t) \} = \{ s,t \}^{\omega}$ and 
$\{ \delta(r), \delta(t) \} = \{ r,t \}^{\pi}$. 
Furthermore, we have $o(\delta(r)\delta(s)) = o(rsrs) =5$.
By the universal property
of Coxeter systems, it follows that $\delta$ extends to an endomorphism
$\alpha$
of $W$. Since $\delta(S)$ generates $W$, $\alpha$ is surjective and hence
an automorphism because $W$ is finite. This finishes (2) and shows in
particular that $\delta(S)$ is a Coxeter generating set. Assertion
(3) is now a consequence of the information collected so far.
\end{proof}

\begin{cor}
Let $(W,S)$, $\omega$, $\pi$ and $\delta$ be as in the previous lemma
and set $c:= rsrs$, $\omega_1 := c\omega$, $\pi_1:=c\pi$ and $\delta_1 := \gamma_c \circ \delta$,
where $\gamma_c$ is the inner automorphism $w \mapsto cwc^{-1}$ of $W$.
Then we have the following:
                                                                                                                             
\begin{itemize}
\item[(1)] $\omega_1 = rsrtsrst$, $\pi_1 = rsrsrts$.
\item[(2)] $\omega_1 s \omega_1^{-1} = srs$, $\omega_1 t \omega_1^{-1} = \pi_1 r \pi_1^{-1}$ and $\pi_1 t \pi_1^{-1}=r$.
\item[(3)] There is an automorphism $\alpha_1$ of $W$ which extends $\delta_1$.
\item[(4)] $\delta_1$ is an $(s,r,rsr)$-deformation of $(W,S)$.
\end{itemize}

\end{cor}
 
\begin{proof} Assertions (1) and (2) are straightforward calculations.
Since $\gamma_c$ is a reflection-preserving automorphism of $W$, Assertions (3) and (4) follow
from Assertions (2) and (3) of the previous lemma, respectively.
\end{proof}

\begin{cor} \label{H3basiccor}
Let $(W,S)$ be a Coxeter system of type $H_3$ where $S = \{ r,s,t \}$
and $o(rs)=5$, $o(rt)=3$. Set $\omega := srstrsrt$, $\pi := srsrstr$
and define $\delta\co S \to W$ by
$\delta(r):= rsr,\delta(s):=s$ and $\delta(t):=\omega t \omega^{-1}$.
Then we have the following:
                                                                                                                             
\begin{itemize}
\item[(1)] $\omega r \omega^{-1} = rsr$, $\omega t \omega^{-1} = \pi s \pi^{-1}$ and $\pi t \pi^{-1}=s$.
\item[(2)] There is an automorphism $\alpha$ of $W$ which extends $\delta$.
\item[(3)] $\delta$ is an $(r,s,srs)$-deformation of $(W,S)$.
\end{itemize}
                                                                                                                             
\end{cor}

\begin{proof} This follows  by exchanging the roles of $r$ and $s$ in the previous corollary.
\end{proof}

\noindent
{\bf Remark:} 
Corollary \ref{H3basiccor} is obtained from Lemma \ref{H3basiclem}
by conjugating by $rsrs$ and then relabelling. We refer to this technique again
in Subsection \ref{SSrelabel} whithout giving further details.

\subsection{Coxeter systems of type $H_4$}\label{omega12subsection}

Throughout this subsection, $(W,S)$ is a Coxeter system of type $H_4$, where $S = \{ r,s,t,u \}$ and $o(rs)=5$, $o(st)=3$. 
Set $J:= \{ r,s \}$, $\omega_1 := rsturstrsrstusrstrs$, $\omega_2 := tsrsrutsrsrtsrsutsrsr$, $\omega_3 := srsrutsrsrtsrsutsrsrtsr$, $\omega := rsrsr\omega_2$, 
$\pi := \omega\omega_1utu$, $\tau := trs\omega_3\omega^{-1}$ and define $\delta\co S \to W$ by
$\delta(r):= rsr$, $\delta(s):=s$, $\delta(t):=\omega t \omega^{-1}$
and $\delta(u)=u$.

\begin{lemma} \label{H4basiclem}
We have the following:
\begin{itemize}
\item[a)] $\pi r \pi^{-1}=rsr$,  
$\omega s \omega^{-1}=s$, $\omega t \omega^{-1} = \pi t \pi^{-1}$ and 
$\omega u \omega^{-1} = u = \pi u \pi^{-1}$.
                                                                               
\item[b)] $\{ \delta(r), \delta(t) \} = \{ r,t \}^{\pi}$,
$\{ \delta(r), \delta(u) \} = \{ r,u \}^{srs}$,
$\{ \delta(s), \delta(t) \} = \{ s,t \}^{\omega}$, \\
$\{ \delta(s), \delta(u) \} = \{ s,u \}^{1_W}$ and
$\{ \delta(t), \delta(u) \} = \{ t,u \}^{\omega}$.

\item[c)] $\tau rsr \tau^{-1}=rsr$, $\tau s \tau^{-1}=s$ and $\tau \omega t \omega^{-1} \tau^{-1}= (tsrtst)t(tsrtst)^{-1}$.

\end{itemize}
\end{lemma}

\begin{proof} The relations in a) and c) are easily deduced from relations given in Franzsen and Howlett \cite[p.333]{FH2}, and b) is an immediate consequence
of a).
\end{proof}

\noindent
{\bf Note:} The relations for $\tau$ will only be needed in Section \ref{tausection}.

\begin{prop} \label{H4basicdefo}
$\delta$ is an $(r,s,srs)$-deformation of $S$ which extends
to an automorphism of $W$.
\end{prop}

\begin{proof} Clearly, $rsr = (srs)r(srs)$ and $s$ generate
$\langle J \rangle$ and $\delta(S)$ generates $W$. 
Setting ${\bf K} := \{ S \}$, it follows that
$\delta$ has Properties a), b) and c) of Proposition \ref{defcritprop}, while Property d) is a consequence of the previous lemma. This proves the claim.
\end{proof}

\subsection{Conditions (TWa)-(TWt)}

Throughout this subsection, $(W,S)$ is a Coxeter system
and $K$ is a subset of $S$ of type $H_k$, where
$k \in \{ 3,4 \}$ and where $r,s,t \in K$ are such
that $o(rs)=5$ and $o(st)=3$; if $k=4$, the unique element
in $K \setminus \{ r,s,t \}$ is denoted by $u$.
Furthermore, we put $J := \{ r,s \}$ and $\omega:=tsrtst$ if $k=3$, $\omega:=rsrsr\omega_2$ if $k=4$,
$\pi:= trs$ if $k=3$ and $\pi:= rsrsr\omega_2\omega_1utu$ if $k=4$, where $\omega_1$ and $\omega_2$ are as in Subsection \ref{omega12subsection}.

\begin{definition} Let $a \in J$.
We say that $K$ is an \emph{$a$-special}  
subset of $S$ or that $K$ is \emph{$a$-special} in
$S$ if the following two Conditions (TWa) and (TWt) hold.
                                                                                
\begin{itemize}
\item[(TWa)] For all $x \in S \setminus K$ we have $o(xa) \in
\{ 2,\infty \}$, and if $o(xa) =2$ then $x \in J^{\perp}$.
\item[(TWt)] If $y \in J^{\perp} \setminus K$ is such that
$o(xy) < \infty$ for some $x \in J^{\infty} \cup \{ t \}$, then $y\in K^{\perp}$.
\end{itemize}
\end{definition}

\begin{lemma} \label{Kspeciallemma}
Let $a \in J$ be such that $K$ is $a$-special in $S$. Then we have
the following.

\begin{itemize}
\item[a)] $\{ K, J^{\infty}, J^{\perp} \setminus K\}$ is a partition
of $S$; if $k=3$ then $K \cap J^{\perp} = \emptyset$ and
if $k=4$ then  $K \cap J^{\perp} = \{ u \}$.
\item[b)] If $y \in J^{\perp}$ is such that
$o(xy) < \infty$ for some $x \in J^{\infty} \cup \{ t \}$, then
$y$ commutes with $\omega$ and with $\pi$.
\end{itemize}
\end{lemma}

\begin{proof} Part a) is immediate and Part b) is a consequence of (TWt) and Lemma \ref{H4basiclem} a).
\end{proof}

\subsection{Angle-deformations for $a$-special subsets of $S$}

We adopt the hypotheses of the previous subsection.
Furthermore, we assume that $a \in J$ is such that
$K$ is $a$-special in $S$.
\\
We define the mapping $\delta\co S \to W$ as follows.
We put $\delta(r) := rsr$, $\delta(y) := y$ 
for all $y \in \{ s \} \cup J^{\perp}$ and $\delta(t) := \omega t \omega^{-1}$.
Let $x \in J^{\infty}$. Then
we put $\delta(x):= \omega x \omega^{-1}$ if $a=r$ and
 $\delta(x):= \pi x \pi^{-1}$ if $a=s$.

\begin{lemma} \label{lemma11.9}
The mapping $\delta$ has the following properties.
                                                                                
\begin{itemize}
\item[a)] $\delta(r) = (srs)r(srs)$ and $\delta(s) = s$ generate
$\langle J \rangle$;
\item[b)] $\delta(S)$ generates $W$;
\item[c)] $\delta \mid_K$ is an $(r,s,srs)$-deformation
of $K$ which extends to an automorphism of $\langle K \rangle$;
\item[d)] for each $x \in S$, there exists an element 
$w_x \in \langle K \rangle$ such that $\delta(x) = w_x x w_x^{-1}$.
\end{itemize}
\end{lemma}

\begin{proof} Assertion a) is obvious. Assertions b) and d)
are immediate consequences of the definition of $\delta$.
Finally, Assertion c) is a consequence of Lemma \ref{H3basiclem}
if $k=3$ and of Proposition \ref{H4basicdefo} if $k=4$.
\end{proof}

\begin{lemma} \label{ordconjlemma}
Let $E$ be an edge of $S$ different from $J$. If
$k=3$ and $a=s$, suppose in addition that
$E$ is not of the form $\{ z,x \}$ with $z \in \{ r,t \}$ and
$x \in J^{\infty}$. Then there exists an element $w_E \in W$
with $\delta(E) = E^{w_E}$.
\end{lemma}

\begin{proof} Let $E = \{ x,y \}$ be such an edge of $S$.
\\
If $E$ is contained in $K$, the assertion follows from Lemma
\ref{H3basiclem} for $k=3$ and Lemma \ref{H4basiclem} for $k=4$.
\\
If $E$ is contained in $J^{\infty}$, then we may choose
$w_E = \omega$ if $a=r$ and $w_E = \pi$ if $a=s$.
\\
If $E$ is contained in $\{ s \} \cup J^{\perp}$, we may choose
$w_E = 1_W$.
\\
If $E$ is contained in $\{ r \} \cup J^{\perp}$, we may choose
$w_E = srs$.

Suppose $E$ is contained in $\{ t \} \cup J^{\perp}$.
As the case $E \subseteq J^{\perp}$ is already
covered by the above, we may assume that $E = \{ t,y \}$
for some $y \in J^{\perp}$.  
Since $o(yt)$ is finite, it follows from Lemma \ref{Kspeciallemma} b)
that $y$ commutes with $\omega$. Hence we may choose
$w_E = \omega$.

Suppose now that $x \in J^{\infty}$ and $y \in J^{\perp}$.
Again by Lemma \ref{Kspeciallemma} b), we know that
$y$ commutes with $\omega$ and with $\pi$. Hence, we
may choose $w_E =  \omega$ if $a=r$ and $w_E = \pi$ if $a=s$. 

Up to renaming the elements of $E$, we are now left 
with the case where $x \in \{ r,s,t \}$ and $y \in J^{\infty}$.

Suppose first that $a=r$. Then the case $x=r$ is not possible
and hence $E$ is contained in $\{ s,t \} \cup J^{\infty}$.
As $s$ commutes with $\omega$ (by Lemma \ref{H4basiclem} a)), we may thus choose $w_E= \omega$.

Suppose now that $a=s$. Then the case $x=s$ is not possible and
by hypothesis, we only have to consider the case $k=4$. In view of the relations given in Lemma \ref{H4basiclem} b), we may choose $w_E = \pi$ in this case, and we are done.
\end{proof}

\begin{prop} \label{defexiprop1}
If $(a,k) \neq (s,3)$, then $\delta$ is an $(r,s,srs)$-deformation
of $S$ which is the restriction of an automorphim of $W$.
\end{prop}

\begin{proof} Setting ${\bf K} = \{ K \}$ in Proposition
\ref{defcritprop}, the two previous lemmas show that
$\delta$ has the required properties and we are done.
\end{proof}

\begin{lemma} \label{switchlemma} 
Suppose $(a,k) = (s,3)$ and let $x \in J^{\infty}$.
Then $\delta(\{ r,x  \}) = \{ t,x \}^{\pi}$ and
 $\delta(\{ t,x  \}) = \{ r,x \}^{\pi}$.
\end{lemma}

\begin{proof} This is an immediate consequence of
the relations given in Lemma \ref{H3basiclem} and the
definition of $\delta$.
\end{proof}

\subsection{$K$-Mirrors}
                                                                                                                             
Throughout this subsection, let $(W,S)$ be a Coxeter system and let
$K = \{ r,s,t \} \subseteq S$ be of type $H_3$ such that
$o(rs) = 5$ and $o(st) = 3$. 
                                                                                                                             
                                                                                            \begin{definition} The \emph{$K$-mirror} of $(W,S)$ is the Coxeter system
$(\bar{W},\bar{S})$ with the property that there
exists a bijection $x \mapsto \bar{x}$ from $S$ onto $\bar{S}$
such that
$o(\bar{r}\thinspace\bar{x}) = o(tx)$ and $o(\bar{t}\thinspace\bar{x}) = o(rx)$ if $x \in J^{\infty}$,
and $o( \bar{x}\thinspace\bar{y}) = o(xy)$ in the remaining cases.
\end{definition}

\noindent
{\bf Remark 1:} Let $(\bar{W},\bar{S})$ be the $K$-mirror of $(W,S)$
and for each $X \subseteq S$, put $\bar{X} := \{ \bar{x} \mid x \in X \}$.
Then $\bar{K}$ is a subset of $\bar{S}$ of type $H_3$ and $(W,S)$ is
the $\bar{K}$-mirror of $(\bar{W},\bar{S})$.

\noindent
{\bf Remark 2:} Let $(\bar{W},\bar{S})$ be the $K$-mirror of $(W,S)$.
Then we have an obvious bijection between the edges of 
$S$ and the edges of $\bar{S}$ which we will call the \emph{canonical
bijection} and which will be denoted by $\theta$.

\noindent
The following lemma is obvious.
\begin{lemma}  \label{mirrorlemma}
Let $(\bar{W},\bar{S})$ be the $K$-mirror of $(W,S)$.
Then $K$ is $s$-special in $S$ if and only if $\bar{K}$ is $\bar{s}$-special in $\bar{S}$. 
\end{lemma}

\subsection{The case $(a,k) = (s,3)$}

Throughout this subsection, let $(W,S)$ be a Coxeter system and let
$K = \{ r,s,t \} \subseteq S$ be of type $H_3$ such that
$o(rs) = 5$ and $o(st) = 3$. We put
$\omega := tsrtst$, $\pi := trs$ and $J:= \{ r,s \}$. Moreover,
$(\bar{W},\bar{S})$ denotes the $K$-mirror of $(W,S)$.
We assume furthermore that $K$ is $s$-special in $S$. Note that this implies
that $\bar{K}$ is $\bar{s}$-special in $\bar{S}$.

We define the mapping $\delta\co S \to W$ by $\delta(x) := x$ if $x \in J^{\perp} \cup \{ s \}$,
$\delta(x) := \pi x \pi^{-1}$ if $x \in J^{\infty}$, $\delta(r):=rsr$ and
$\delta(t) := \omega t \omega^{-1}$.

\begin{lemma} \label{sTWlem2}
Let $\{ x,y \}$ be an edge of $S$. Then $o(\delta(x) \delta(y)) = o(\bar{x}\thinspace\bar{y})$.
\end{lemma}

\begin{proof} 
This is a consequence of Lemmas \ref{ordconjlemma} 
and \ref{switchlemma}.
\end{proof}

\begin{lemma} \label{sTWlem3}
$\delta(S)$ is a Coxeter generating set of $W$. Moreover, there exists a bijection $\Delta$ from the set of edges of $S$ onto the set of edges of $\delta(S)$ such that
$\Delta(J) = \{  rsr,s  \}$ and such that
for each edge $E \neq J$ of $S$, there exists $w_E \in W$
with $\Delta(E) = E^{w_E}$. 
\end{lemma}

\begin{proof} By the universal property of $(\bar{W},\bar{S})$
and Lemmas \ref{sTWlem2} and \ref{lemma11.9}, there is an epimorphism
$\bar{\beta}\co \bar{W} \to W: \bar{x}\mapsto \delta(x)$ with the following
properties:
\begin{itemize}
\item[a)] $\bar{\beta} \mid_{\langle \bar{K} \rangle}$ is an
isomorphism from $\langle \bar{K} \rangle$ onto 
$\langle K \rangle$.
\item[b)] For each $x \in S$, there exists an element $w_x \in \langle K \rangle$
such that $\bar{\beta}(\bar{x}) = w_x x w_x^{-1}$.
\end{itemize}

By Lemma \ref{mirrorlemma}, $\bar{K}$ is 
 $\bar{s}$-special in $\bar{S}$. Hence, by defining
$\bar{\omega},\bar{\pi} \in \bar{W}$ and 
$\bar{\delta}\co \bar{S} \to \bar{W}$ for $(\bar{W},\bar{S})$,
we obtain also an epimorphism $\beta\co W \to \bar{W}$
with the following properties:
\begin{itemize}
\item[a)] $\beta \mid_{\langle K \rangle}$ is an
isomorphism from $\langle K \rangle$ onto
$\langle \bar{K} \rangle$.
\item[b)] For each $x \in S$, there exists an element $\bar{w}_x \in \langle \bar{K} \rangle$
such that $\beta(x) = \bar{w}_x \bar{x}\thinspace \bar{w}_x^{-1}$.
\end{itemize}

We put $\alpha:= \bar{\beta} \circ \beta$ and for each $x \in S$,
we set $v_x := \bar{\beta}(\bar{w}_x) w_x$. Then $\alpha\co W \to W$
is an epimorphism with the following properties:
\begin{itemize}
\item[a)] $\alpha \mid_{\langle K \rangle}$ is an
automorphism $\langle K \rangle$.
\item[b)] For each $x \in S$, we have 
$v_x \in \langle K \rangle$ and $\alpha(x) = v_x x v_x^{-1}$. 
\end{itemize}
                               
Now, it follows from Proposition \ref{autoprop} (with ${\bf K}=\{K\}$) that
$\alpha$ is an automorphism of $W$. In particular,
$\bar{\beta}$ is an isomorphism. As $\delta(S) = \bar{\beta}(\bar{S})$,
the set $\delta(S)$ is a Coxeter generating set of $W$.          

It remains to find an appropriate $\Delta$.
As $\bar{\beta}$ is an isomorphism, we have a canonical bijection $\Delta_1$ from
the set of edges of $\bar{S}$ onto the set of edges of $\delta(S)$. Let $\theta$ be the canonical bijection from the set of edges of $S$ onto the set of edges of $\bar{S}$.
It is then readily verified, using Lemma \ref{ordconjlemma} 
and \ref{switchlemma}, that $\Delta := \Delta_1 \circ \theta$ is the required bijection. This finishes the proof of the lemma. 
\end{proof}

\begin{prop} \label{sTWprop}
$\delta$ is an $(r,s,srs)$-deformation of $(W,S)$.
\end{prop}

\begin{proof} This is a consequence of the two previous lemmas.
\end{proof}

\subsection{The relabeled version} \label{SSrelabel}

Throughout this subsection, $(W,S)$ is a Coxeter system
and $K$ is a subset of $S$ of type $H_k$, where
$k \in \{ 3,4 \}$ and where $r,s,t \in K$ are such
that $o(rs)=5$ and $o(rt)=3$; if $k=4$, the unique element
in $K \setminus \{ r,s,t \}$ is denoted by $u$.
Define $\overline{\omega}_i$ for $i\in\{1,2\}$ by exchanging r and s in the expression of $\omega_i$ given in Subsection \ref{omega12subsection}, where $t$ and $u$ are as above. Also, let $c:=rsrs$ and $\overline{c}:=srsr$.
We put $J := \{ r,s \}$, $\omega:=\overline{c}trstrt = srstrsrt$ if $k=3$, $\omega:=\overline{c}srsrs\overline{\omega}_2 = r\overline{\omega}_2$ if $k=4$,
$\pi:= \overline{c}tsr = srsrstr$ if $k=3$ and $\pi:= \overline{c}srsrs\overline{\omega}_2\overline{\omega}_1utu = r\overline{\omega}_2\overline{\omega}_1utu$ if $k=4$.
We assume that $a \in J$ is such that $K$ is $a$-special in $S$
and we define $\delta\co S \to W$ as follows.
We put $\delta(r) := rsr$, $\delta(y) := y$
for all $y \in \{ s \} \cup J^{\perp}$ and $\delta(t) := \omega t \omega^{-1}$.
Let $x \in J^{\infty}$. Then
we put $\delta(x):= \omega x \omega^{-1}$ if $a=s$ and
 $\delta(x):= \pi x \pi^{-1}$ if $a=r$.
                                                                                                                             
The following proposition is obtained from Propositions 
\ref{defexiprop1} and \ref{sTWprop}
by relabelling.

\begin{prop} \label{relabelprop}
The mapping $\delta$ is an $(r,s,srs)$-deformation of $(W,S)$.
Moreover, if $(k,a) \neq (3,r)$, it is the restriction of an automorphism
of $W$.
\end{prop}

\section{$\Delta$-edges}\label{figDE}

\subsection{Some particular diagrams}

Throughout this subsection, we put $\lambda := 2\cos(\pi/5)$.

\begin{figure}
 \begin{minipage}[b]{.46\linewidth}
  \centering\epsfig{trim = 0mm 60mm 0mm 0mm, clip, height=6cm, width=4cm, figure=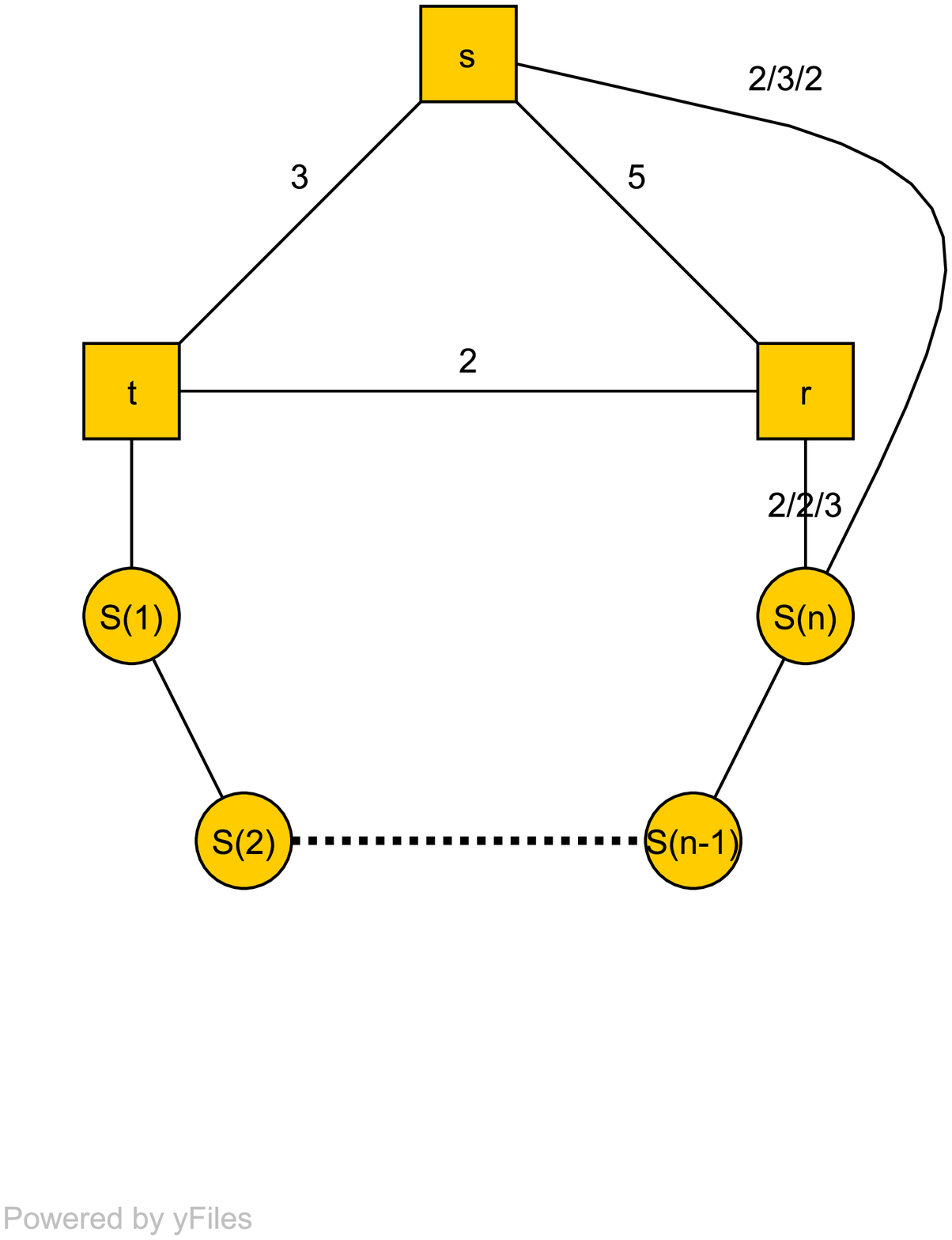,width=\linewidth}
  \caption{{\bf (DE3)}: $n\geq 2$ and $\overline{X_1t}=\infty$ and $X_n\subseteq \{r,s\}^{\infty}.$ \label{figDE3}}
 \end{minipage} \hfill
 \begin{minipage}[b]{.46\linewidth}
  \centering\epsfig{trim = 0mm 60mm 0mm 0mm, clip, height=6cm, width=4cm, figure=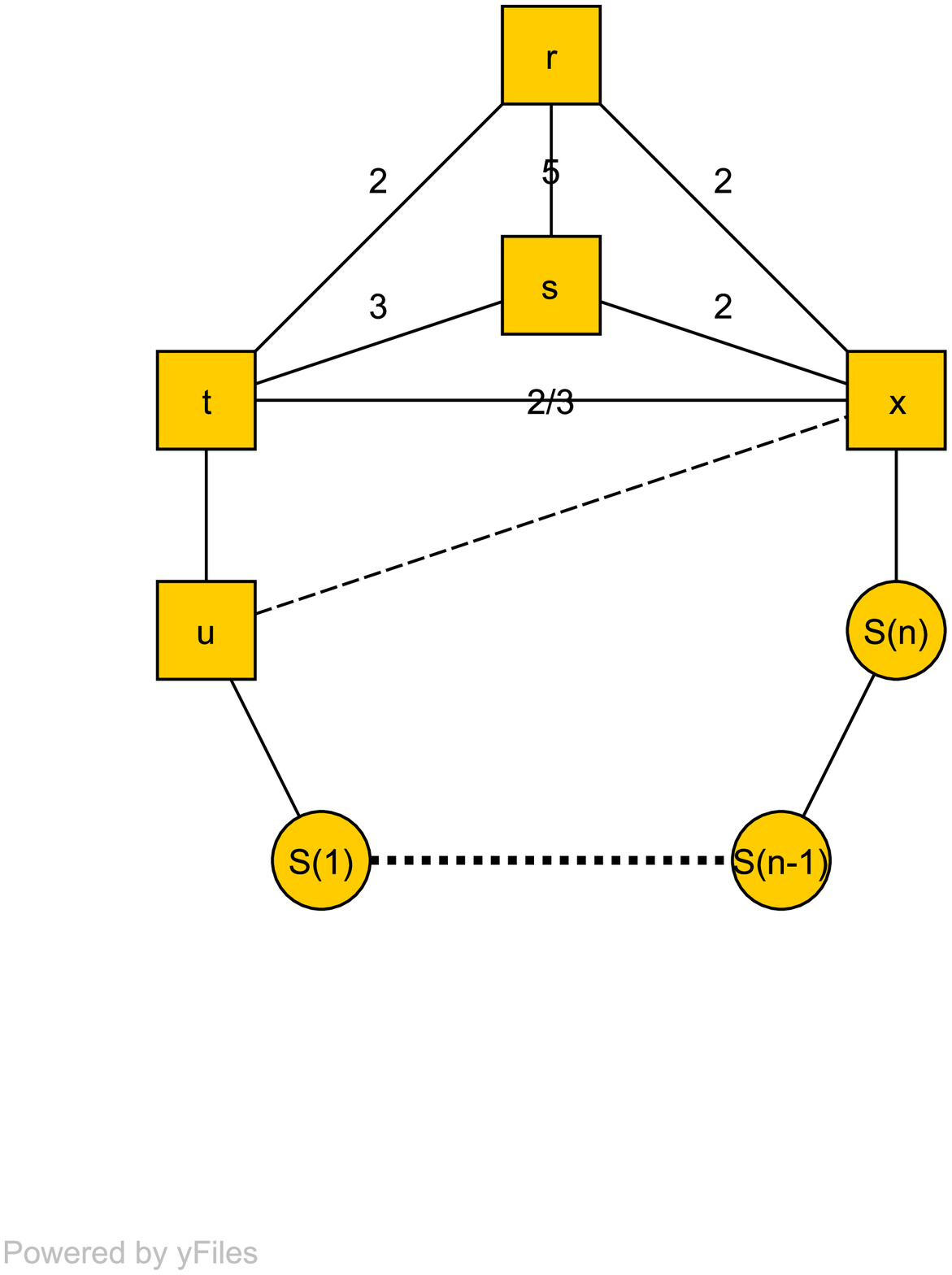,width=\linewidth}
  \caption{{\bf (DE4)}: $n\geq 2$ and $\overline{X_1t}=\infty$ and $X_n\subseteq \{r,s\}^{\infty}.$ \label{figDE4}}
 \end{minipage}
\end{figure}

\noindent
Let $(W_1,R_1)$ be a Coxeter system whose diagram
is as in Figure \ref{figDE3} and
let $(W_2,R_2)$ be a Coxeter
system whose diagram is as in Figure \ref{figDE4}.
Hence, we have $R_1 = \{ r,s,t \} \cup X$ and
$R_2 = \{ r,s,t,u,x\} \cup X$ where
$X = \{ S(i) \mid 1 \leq i \leq n \}$.
                                                                                
\noindent
For $k=1,2$, we
consider the geometric representation of $(W_k,R_k)$
and its root system $\Phi_k$;
in particular we identify $W_k$ with its image in $O(V_k,b_k)$.

\noindent
We put $\alpha_1 := rs(e_r)=  \lambda e_r + \lambda e_s$,
$\Pi_1 := \{ \alpha_1,e_t \} \cup \{ e_{S(i)} \mid 1 \leq i \leq n \}$,
$S_1:=\{ \rho_{\alpha} \mid \alpha \in \Pi_1 \}$ and $\omega_1 := rst$.

\noindent
We put $\alpha_2 := srstrs(e_r)=  (\lambda+1) e_r + 2\lambda e_s + \lambda e_t$, $\Pi_2 := \{ \alpha_2,e_u,e_x \} \cup \{ e_{S(i)} \mid 1 \leq i \leq n \}$,
$S_2:=\{ \rho_{\alpha} \mid \alpha \in \Pi_2 \}$ and $\omega_2 := srstrsut$.

\noindent
The following facts are easily verified for $k=1,2$:
                                                                                
\begin{itemize}
\item[a)] $\rho_{\alpha_k} = \omega_k r \omega_k^{-1}$.
\item[b)] $\Pi_k$ is a root subbase of $\Phi_k$;
in particular $S_k$ is a fundamental set.
\item[c)] $\Gamma(S_k)$ is a chordfree circuit.
\item[d)] $\omega_1 s \omega_1^{-1} = t$
and  $\omega_2 s \omega_2^{-1} = u$.
\end{itemize}

\subsection{Coxeter systems containing some particular subsystems}
Throughout this subsection, $(W,R)$ is a Coxeter system
and $W$ is identified with its image in $O(V,b)$ via its geometric representation.

\begin{prop} \label{Figprop}
For $k=1,2$, let $R_k \subseteq R^W$ be a fundamental
set of reflections and put $W_k := \langle R_k \rangle$.
Suppose that $(W_k,R_k)$ is a Coxeter system whose diagram
is as in Figure \ref{figDE3} if $k=1$ and as in Figure
\ref{figDE4} if $k=2$.
Then $\{ r,s \}$ is sharp-angled.
\end{prop}
                                                                                
\begin{proof} For $k=1,2$, we define $\omega_k \in W_k$
as in the previous subsection.
We put $S_1 := (R_1 \setminus \{ r,s \}) \cup \{ \omega_1 r \omega_1^{-1} \}$
and $S_2 := (R_2 \setminus \{ r,s,t \}) \cup \{ \omega_2 r \omega_2^{-1} \}$.
By the considerations above, we know that
the set $S_k$  is
a fundamental set of reflections. Moreover, $\Gamma(S_k)$
is a chordfree circuit. By Theorem \ref{reflrigidthm}, it
follows that $S_k$ is sharp-angled.
Hence  $\{ \omega_1 r \omega_1^{-1},t \}$ and
$\{ \omega_2 r \omega_2^{-1},u \}$ are sharp-angled.
As $\omega_k$ is an element
of $W$ which conjugates $\{ r,s \}$ onto $\{ \omega_1 r \omega_1^{-1},t \}$
for $k=1$, and onto $\{ \omega_k r \omega_k^{-1},u \}$ for $k=2$, 
it follows that $\{ r,s \}$ is sharp-angled as well.
\end{proof}

\begin{cor}\label{figcor}
Let $S \subseteq R^W$ be a fundamental set of 
reflections and let $J = \{ r,s \}$ be an edge of $S$ such that $o(rs)=5$, and which is not sharp-angled. Then
there is no subset $K$ of $S$ as in Figures
\ref{figDE3} or \ref{figDE4}.
\end{cor}

\subsection{Definition of $\Delta$-edges}

\begin{definition}
Let $W$ be a group and $S$ a subset of involutions of $W$. Let $J = \{ r,s \}$
be an edge of $S$. We call $J$ a \emph{$\Delta$-edge} of $S$ if there is no subset $K$ of $S$ containing $J$ having one of the following properties:
                                                                                
\begin{itemize}
\item[(DE1)]
$\Gamma(K)$ is non-spherical, 2-spherical and irreducible.
\item[(DE2)]
$\Gamma(K)$ is a chordfree circuit of length at least 4.
\item[(DE3)]
$\Gamma(K)$ is a diagram as shown in figure \ref{figDE3}.
\item[(DE4)]
$\Gamma(K)$ is a diagram as shown in figure \ref{figDE4}.
\end{itemize}
\end{definition}

\noindent
{\bf Remark:} Note that if $o(rs) \neq 5$, then $J$ is a $\Delta$-edge if and only if (DE1) and (DE2) hold; if $o(rs)=5$, the same remains true if there is
no subset of type $H_3$ containing $J$.

\noindent The definition of $\Delta$-edges is motivated by the following
proposition, which is a consequence of Theorem \ref{reflrigidthm}
and Corollary \ref{figcor}.

\begin{prop}\label{deftheorem}
Let $(W,R)$ be a Coxeter system, let $S \subseteq R^W$ be a fundamental
set of reflections and suppose that $J$ is an edge  of $S$ which is
not sharp-angled with respect to $R$. Then $J$ is a $\Delta$-edge of $S$.
\end{prop}

\section{$\Delta$-edges of type $H_2$}
\label{sec9}

Throughout this section, $(W,S)$ is a Coxeter system
and $J = \{ r,s \} \subseteq S$ is a $\Delta$-edge
of $(W,S)$ with $o(rs)=5$. Moreover, we define several
subsets of $S$ as follows.

\begin{itemize}
\item
$T := \{t\in S \ | \ \textrm{type($\{r,s,t\})=H_3$}\} = T_r \amalg T_s$
\end{itemize}
where $T_r := \{t\in T \ | \ m_{rt}=3\}$ and $T_s := \{t\in T \ | \ m_{st}=3\}$.
\begin{itemize}
\item
For a $J$-component $L$, put $T_L := \{t\in T \ | \ \exists x\in L \ : \ m_{xt}<\infty\}$.
\item
$U := \{u\in S \ | \ \exists t\in T \ \textrm{such that type($\{r,s,t,u\}$)}= H_4\}$.
\item
For $t\in T$, set $U_t := \{u\in U \ | \ \textrm{type($\{r,s,t,u\}$)}= H_4\}$.
\item
For $t\in T$ and $L$ a $J_t$-component, let $U_L := \{u\in U_t \ | \ \exists x\in L \ : \ m_{xu}<\infty\}$.
\item
For $t\in T$ and $u\in U_t$, let $J_t:=J\cup \{t\}$ and $J_{t,u}:= J\cup \{t\}\cup \{u\}$.
\item
For $u\in U$, set $T_u := \{t\in T \ | \ \textrm{type($\{r,s,t,u\}$)}= H_4\}$.
\item
$T^3 := \{t\in T \ | \ U_t = \emptyset\}$.
\item
$T^4 := T \setminus T^3$.
\item
For $a \in J$ and $k \in \{ 3,4 \}$, put
$T_a^k := T_a \cap T^k$.
\end{itemize}

\subsection{Some preliminary observations}

\begin{lemma}\label{Jflexlemma}
$J$ is flexible.
\end{lemma}

\begin{proof} This is Lemma \ref{flexedgelemma}.
\end{proof}

\begin{lemma} \label{noedgesinT}
There are no edges in $T$ and for each $t \in T$,
there are no edges in $U_t$.
\end{lemma}

\begin{proof} This follows from (DE1).
\end{proof}

\subsection{Flexibility of $J_t$ and consequences}\label{t(L)subsection}

\begin{prop} \label{Jtflexprop}
For all $t \in T$, the set $J_t$ is flexible.
\end{prop}

\begin{figure}
\begin{center}
\includegraphics[trim = 0mm 55mm 0mm 0mm, clip, height=5cm]{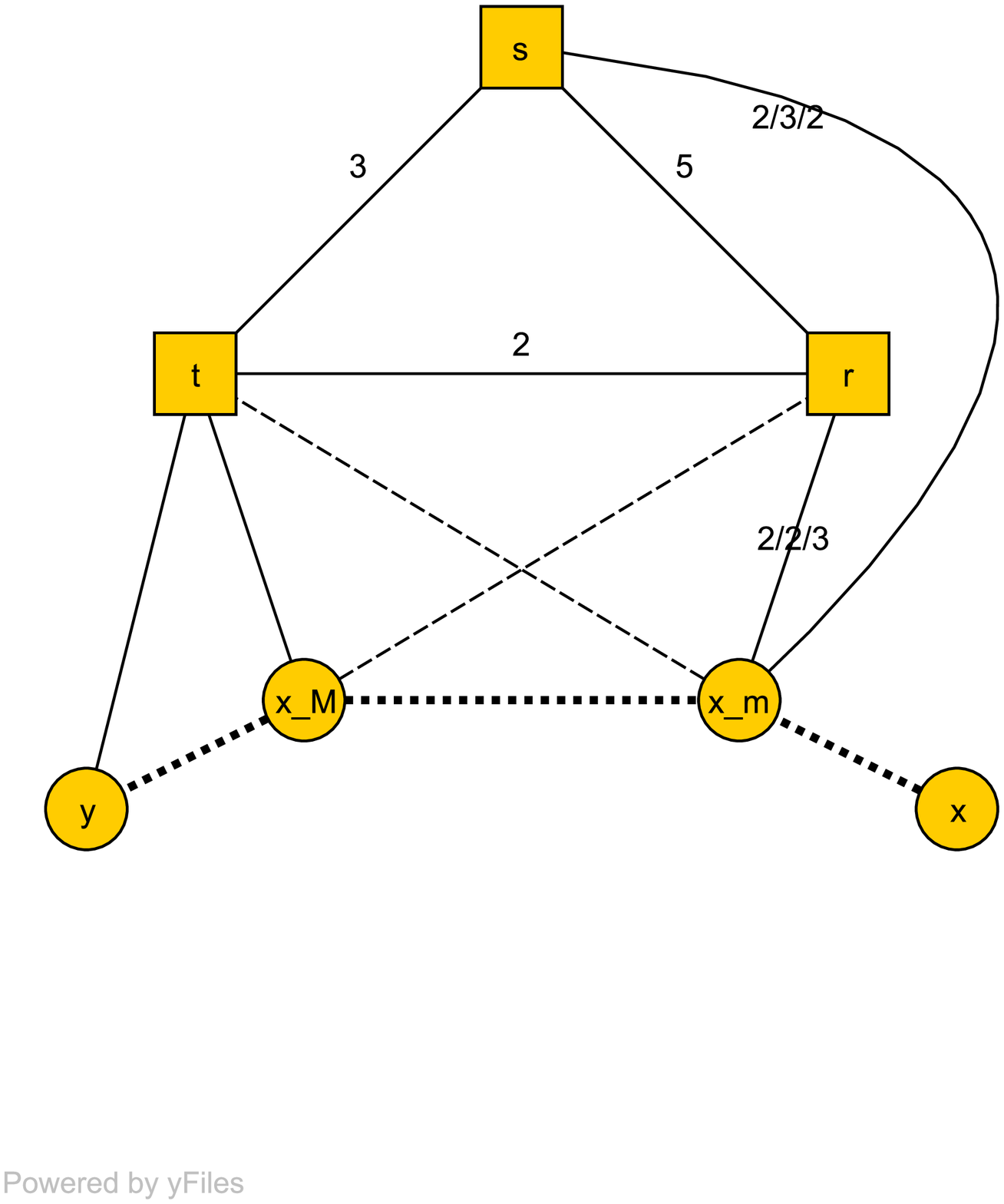}
\caption{Contradicts (DE3).}
\label{graph2_Jtflexible}
\end{center}
\end{figure}

\begin{proof}
Let $t\in T$ and let $L$ be a $J_t$-component. If $L$ is also a $J$-component, then $L$ is flexible by Lemma \ref{Jflexlemma} and we are done. So, we may assume there exists an $x\in L$ such that $x\in J^{\fin}$ (thus $m_{xt}=\infty$). Suppose by contradiction there exists $y\in L$ such that $m_{yt}<\infty$. Then $m_{yr}=\infty$ or $m_{ys}=\infty$.
\\
Let $x=x_0, x_1, \dots, x_k = y$ be a minimal path in $L$ joining $x$ to $y$. Define
$$M := min\{i \ | \ 0< i \leq k; \ m_{x_it}<\infty \}$$
and
$$m := max\{i \ | \ 0\leq i < M; \ x_i\in J^{\fin}\}.$$
Then the subpath $x_m, x_{m+1}, \dots , x_M$ from $x_m$ to $x_M$ is still minimal, hence chordfree, and possesses the following properties:
\begin{itemize}
\item[(1)]
$(m_{x_mr},m_{x_ms}) \in \{(2,2), (2,3), (3,2)\}$ (by (DE1));
\item[(2)]
$m_{x_it}=\infty$ for all $i$ such that $m\leq i<M$ (by definition of M);
\item[(3)]
$x_i\in J^{\infty}$ for all $i$ such that $m<i \leq M$ (by definition of m).
\end{itemize}
Moreover, $m_{x_Mt}<\infty$. In conclusion, we obtain a subgraph $\{r, x_m, x_{m+1}, \dots, x_M, t, s\}$ as pictured in figure \ref{graph2_Jtflexible}, contradicting (DE3).
\end{proof}
                                                                                
\begin{cor} \label{Jtflexcor1}
Let $t\in T$ and let $L$ be a $J$-component such that
there exists $z \in L$ with $o(zt) < \infty$. If $y \in J^{\fin}\setminus \{t\}$
is such that there exists an $x  \in L$ with $o(xy) < \infty$, then
$y \in J_t^{\fin}$.
\end{cor}

\begin{proof}
Let $L'$ be the $J_t$-component containing $L$.
If $o(yt)=\infty$, we get $y \in L'$ because $o(xy) < \infty$.
But then $z$ and $y$ belong to $L'$, contradicting
the fact that $J_t$ is flexible. Hence $o(yt)<\infty$ and
so $y \in J_t^{\fin}$ because $y \in J^{\fin}$ by
assumption.
\end{proof}

\begin{cor} \label{Jtflexcor2}
Let $L$ be a $J$-component, then $|T_L| \leq 1$.
\end{cor}

\begin{proof} This follows from the previous corollary
and Lemma \ref{noedgesinT}.
\end{proof}

\begin{definition} Let $L$ be a $J$-component. If $T_L$
is non-empty, then $t(L)$ denotes its unique element;
if $T_L$ is empty, we put $t(L) := \infty$.
\end{definition}

\subsection{Flexibility of $J_{t,u}$ and consequences}

\begin{prop} \label{Jtuflexprop}
Let $t \in T$ and $u \in U_t$. Then $J_{t,u}$ is flexible.
\end{prop}

\begin{figure}
\begin{center}
\includegraphics[trim = 0mm 60mm 0mm 0mm, clip, height=5cm]{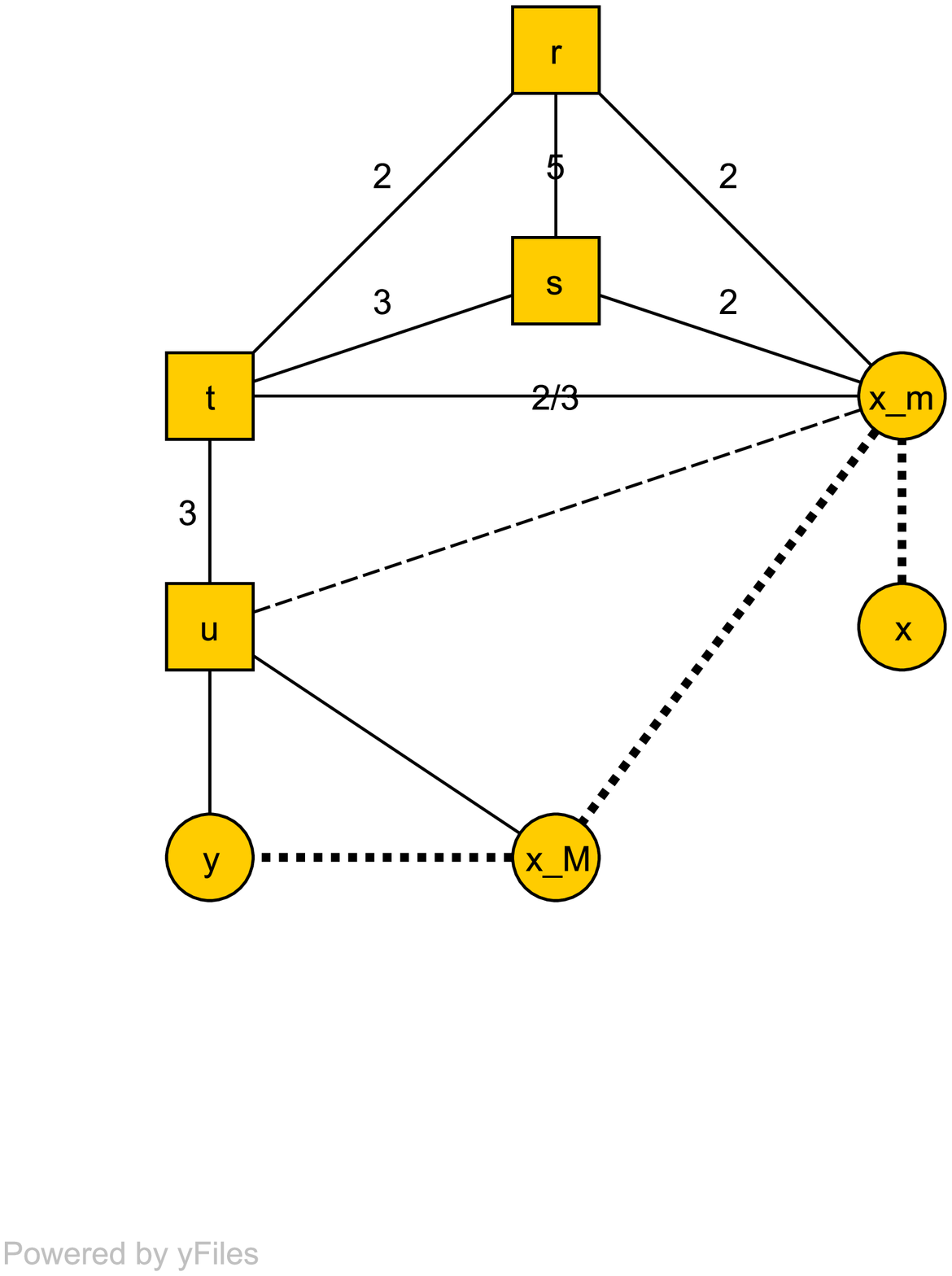}
\caption{Contradicts (DE4).}
\label{graph3_Jtuflexible}
\end{center}
\end{figure}

\begin{proof}
Let $L$ be a $J_{t,u}$-component. If $L$ is also a $J_t$-component, then it is free by Proposition \ref{Jtflexprop} and we are done. So, we may assume there exists
an $x\in L$ such that $x\in J_t^{\fin}$ (thus $m_{xu}=\infty$). Suppose by contradiction there exists $y\in L$ such that $m_{yu}<\infty$. Then $y\in J_t^{\infty}$.
\\
Let $x=x_0, x_1, \dots, x_k = y$ be a minimal path in $L$ joining $x$ to $y$. Define
$$M := min\{i \ | \ 0< i \leq k; \ m_{x_iu}<\infty \}$$
and
$$m := max\{i \ | \ 0\leq i < M; \ x_i\in J_t^{\fin}\}.$$
Then the subpath $x_m, x_{m+1}, \dots , x_M$ from $x_m$ to $x_M$ is still minimal, hence chordfree, and possesses the following properties:
\begin{itemize}
\item[(1)]
$(m_{x_mr}, m_{x_ms}, m_{x_mt}) \in \{(2,2,2), (2,2,3)\}$ (by (DE1));
\item[(2)]
$m_{x_iu}=\infty$ for all $i$ such that $m\leq i<M$ (by definition of M);
\item[(3)]
$x_i\in J_t^{\infty}$ for all $i$ such that $m<i \leq M$ (by definition of m).
\end{itemize}
Moreover, $m_{x_Mu}<\infty$. In conclusion, we obtain a subgraph $\{r, x_m, x_{m+1}, \dots, x_M, u, t, s\}$ as pictured in figure \ref{graph3_Jtuflexible}, contradicting (DE4).
\end{proof}

\begin{cor} \label{H4flexcor1}
Let $t \in T$, $u \in U_t$ and $L$ be a $J_t$-component containing
an element $z$ with $o(zu) < \infty$. Suppose that $y \in J_t^{\fin}$
is such that there exists $x \in L$ with $o(xy) < \infty$.
Then $y \in J_{t,u}^{\fin}\cup\{u\}$; in particular, if $y\neq u$, then $y \in J_{t,u}^{\perp}$ .
\end{cor}

\begin{proof}
Let $L'$ be the $J_{t,u}$-component containing $L$ and suppose $y\neq u$.
If $o(yu)=\infty$, we get $y \in L'$ since $o(xy) < \infty$.
But then $z$ and $y$ belong to $L'$, contradicting
the flexibility of $J_{t,u}$. Hence $o(yu)<\infty$ and
so $y \in J_{t,u}^{\fin}$ because $y \in J_t^{\fin}$ by
assumption. Now, (DE1) implies that $J_{t,u}^{\fin} = J_{t,u}^{\perp}$, so we are done.
\end{proof}

\begin{cor} \label{Jtuflexcor2}
Let $t \in T$ and let $L$ be a $J_t$-component.
Then $|U_L| \leq 1$.
\end{cor}

\begin{proof}
This follows from the previous corollary.
\end{proof}

\begin{definition} Let $t \in T$ and let $L$ be a $J_t$-component. If $U_L$
is non-empty, then $u(L)$ denotes its unique element;
if $U_L$ is empty, we put $u(L) := \infty$.
\end{definition}

\noindent
{\bf Remark:} Let $t \neq t' \in T$. By Lemma \ref{noedgesinT},
we can talk about the `$J_t$-component containing $t'$' as we will
do in the following proposition.

\begin{prop} \label{tneqt'prop} 
 Let $t \neq t' \in T$, let $L$ be the $J_t$-component containing $t'$
and put $K := J_t \cup U_L$.
Then $J_{t'}^{\fin}$ is contained in $K^{\fin} \cup L \cup U_L$.
\end{prop}

\begin{proof} Let $y \in J_{t'}^{\fin}$. Then we have
in particular $o(yt') <  \infty$. Hence, if $o(yt) = \infty$,
we have $y \in L$. Thus we are left with the case
where $o(yt) < \infty$. As $y\in J_{t'}^{\fin}$, we get that $y \in J_t^{\fin}$.
In particular, we are already done if $u(L) = \infty$.

Let us now assume that $U_L \neq \emptyset$ and put $u:=u(L)$.
Then there exists an element $z \in L$ such that
$o(uz) < \infty$ and there exists an element
$x \in L$ (namely $t'$) such that $o(xy)  < \infty$.
As $y\in J_t^{\fin}$, the claim follows from
Corollary \ref{H4flexcor1}.
\end{proof}

\subsection{Tameness}

\begin{figure}
\begin{center}
\includegraphics[trim = 0mm 120mm 0mm 0mm, clip, height=4cm]{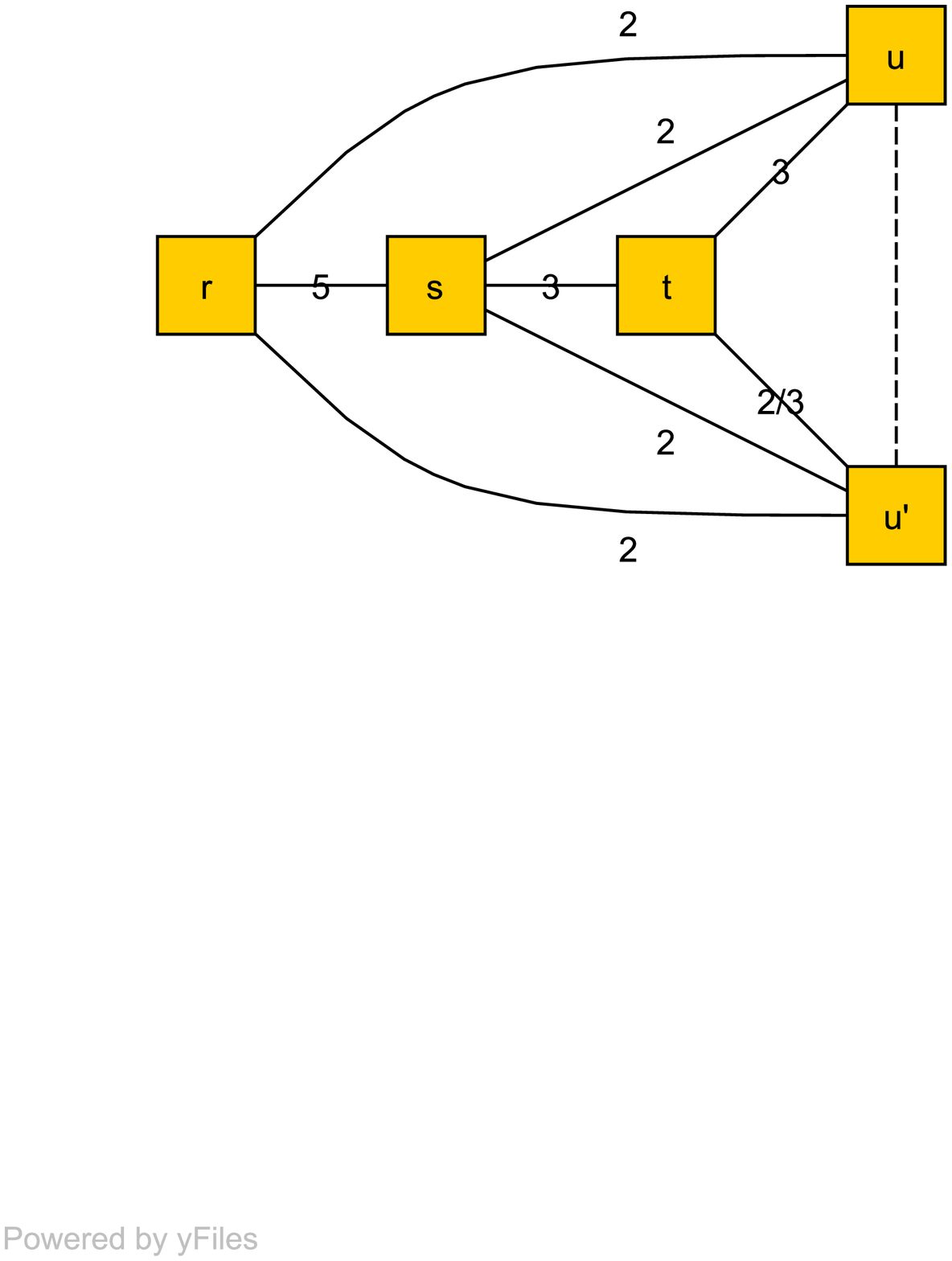}
\caption{Tameness.}
\label{graph5_hyp}
\end{center}
\end{figure}

\begin{definition} Let $t \in T$ and let $K$ be a subset
of $S$ containing $J_t$. Then $t$ is called \emph{tame} in $K$
if there is no subset $K'$ of $K$ containing $J_t$ such that $\Gamma(K')$ is as in Figure \ref{graph5_hyp}. We call $t$ \emph{tame}, if it is tame in $S$. Otherwise, we call it \emph{wild}.
\end{definition}

\noindent
Here are some basic observations. The first two of them
are obvious whereas the third one is a consequence of Lemma \ref{noedgesinT}.

\begin{lemma}
Let $t \in T$ and $K_1 \subseteq K$ be  subsets of $S$ containing $J_t$.
If $t$ is tame in $K$, then it is tame in $K_1$.
\end{lemma}

\begin{lemma}
If $t \in T^3$, then $t$ is tame.
\end{lemma}

\begin{lemma}\label{Utleq1lemma}
If $t$ is tame, then $|U_t| \leq 1$.
\end{lemma}

\noindent
Let $t \in T$ be tame. Then we put $K_t := J_t \cup U_t$.

\begin{lemma} \label{lemma9.13}
Let $t \in T$ be tame. Then $J_t^{\perp} = K_t^{\perp}$
and   $J_t^{\fin}\cup J_t = K_t \cup K_t^{\perp}$.
\end{lemma}

\begin{proof} We start with the first equality which is
trivial if $U_t$ is empty. Suppose $U_t$ is non-empty 
and let $u$ denote its unique element. 
Obviously, we have $K_t^{\perp} \subset J_t^{\perp}$.
Let now $y \in J_t^{\perp}$. If $o(yu) = \infty$, we get
a contradiction to the tameness of $t$ (using (DE1)) and if
$2< o(yu) < \infty$, we get a contradiction to (DE1).
Hence $o(uy) =2$ and the first equality holds.

The second equality follows now from the fact that
$J_t^{\fin} = J_t^{\perp} \cup U_t$ (because of (DE1)),
the definition of $K_t$ and the first equality.
\end{proof}

\begin{lemma} \label{lemma9.14}
Let $t \in T$ be tame, $K:=K_t$, let $L$ be a $J$-component
with $t = t(L)$ and let $a \in J$ be $L$-free.
Then $K$ is an $a$-special subset of $S':= K \cup L \cup J^{\perp}$.
\end{lemma}

\begin{proof} Note first that $S' \setminus K \subseteq L \cup J^{\perp}$. 
Thus, as $a$ is $L$-free, Condition (TWa) is obviously satisfied.

We now show that Condition (TWt) holds as well. Note first
that $J^{\infty} \cap S' = L$. 
Let $y \in  J^{\perp}
\setminus K$ such that $o(yx) < \infty$ for some 
$x \in (J^{\infty} \cap S')\cup \{t\} = L \cup \{ t \}$. We first show that
$y \in  J_t^{\fin}$, which is obvious if $x =t$. Hence
we may assume $x \in L$. 
As $t = t(L)$, there
exists $z \in L$ such that $o(tz) < \infty$.  
Therefore, $y \in  J^{\perp} \subseteq J^{\fin}\setminus \{t\}$ and we can apply Corollary \ref{Jtflexcor1} to see that
$y \in J_t^{\fin}$. 

Now, as $t$ is tame and $y$ is not in $K$, we have
$y \in  J_t^{\perp}$ and we are done if $U_t = \emptyset$.
Suppose $U_t \neq \emptyset$ and let $u$ be the unique 
element of $U_t$. If $o(yu) = \infty$, we get a contradiction
to the tameness of $t$ and if $2 < o(yu) < \infty$, we get
a contradiction to (DE1). Hence $o(yu) = 2$ and $y\in K^{\perp}$ because $K = J_t \cup \{ u \}$ and
$y \in  J_t^{\perp}$.
\end{proof}

\subsection{The degree of a subset containing $J$}\label{WXYZsubsection}

\begin{definition} Let $K$ be a subset of $S$ containing
$J$. The \emph{degree} of $K$ is the number of elements in $K \cap T$
which are wild in $K$. It is denoted by $\deg(K)$.
\end{definition}

\noindent
Here is a preliminary observation.

\begin{lemma}
Let $J \subseteq K_1 \subseteq K \subseteq S$. Then $\deg(K_1) \leq \deg(K)$.
\end{lemma}

\begin{definition} Let $t \in T$. For each
$u \in \widehat{U}_t := U_t \cup \{ \infty \}$, we define the
sets $V_u,W_u,X_u,Y_u$ and $Z_u$ as follows.

\begin{itemize}
\item
$V_{\infty} := J_t$ and $V_u := J_{t,u}$ for $u \in U_t$;
\item
$W_u := V_u \cup V_u^{\perp}$;
\item
$X_u$ is the union of all $J_t$-components $L$ such that
$u(L) = u$;
\item
$Y_u := W_u \cup X_u$;
\item
$Z_u := Y_u \cup Y_{\infty}$.
\end{itemize}
\end{definition}

\begin{lemma}\label{deglemma}
Let $t \in T$ and $u \in \widehat{U}_t$. Then $t$ is tame
in $Y_u$. In particular, if $t$ is wild then
$\deg(Y_u) < \deg(S)$.
\end{lemma}

The following is a consequence of Proposition \ref{tneqt'prop}.

\begin{lemma} \label{lemma922}
Let $t \neq t' \in T$ and $u \in \widehat{U}_t$. 
If $t'$ is contained in $X_u$, then $J_{t'}^{\fin} \subseteq Y_u$.
\end{lemma}

\begin{lemma} \label{lemma923}
Let $u \in U_t$. Then $Y_u \cap Y_{\infty} = J_t\cup J_{t,u}^{\perp}$
and if $E$ is an edge of $Z_u$, then $E \subseteq Y_u$ or 
$E \subseteq Y_{\infty}$.
\end{lemma}

\begin{proof} The first statement follows from the definition
of the sets $Y_u$ and $Y_{\infty}$.
\\
Let $E = \{ x,y \}$ be an edge of $Z_u$ and suppose
that $x \in Y_u$ and $y \in Y_{\infty}$. 

Suppose first that $x \in X_u$.
Then $y$ cannot be in $X_{\infty}$ since in that case $x$ and $y$ would be in
different $J_t$-components. Hence, $y \in J_t \cup J_t^{\perp}$.
If $y \in J_t \cup J_{t,u}^{\perp}$, then $y$ is in $Y_u$
and we are done. Suppose by contradiction that $y \in J_t^{\perp} \setminus  J_{t,u}^{\perp}$.
Then we have $o(yu)=\infty$ by (DE1).
Let $L$ be the $J_t$-component containing $x$.
Then there is an element $z$ in $L$ such that
$o(uz)$ is finite. Let $L'$ be the $J_{t,u}$-component
containing $L$. Then $x,y$ and $z$ are contained in $L'$,
contradicting the flexibility of $J_{t,u}$.

Thus we may assume that $x \in J_{t,u} \cup J_{t,u}^{\perp}$.
If $x \neq u$, we have $x \in Y_{\infty}$ and we are done.
Suppose that $x=u$. Then the case $y \in X_{\infty}$
is not possible, because otherwise we would have
$u = u(L)$ for the $J_t$-component $L$ containing $y$.
Thus we may assume that $y \in J_t \cup J_t^{\perp}$.
By (DE1), we then get $y \in J_t \cup J_{t,u}^{\perp}$ and hence $E \subseteq Y_u$ in this case. 
\end{proof}

\section{Existence of Angle-deformations}\label{tausection}

Throughout this section, $(W,S)$ is a Coxeter system
and $J = \{ r,s \} \subseteq S$ is a $\Delta$-edge 
of $(W,S)$ with $o(rs)=5$.

\noindent
We adopt the notations of the previous section.

\subsection{Conventions for tame elements and standard deformations}
\label{sstameconv}

If $t \in T$ is tame, we fix the following
notations:
\begin{itemize}
\item
By Lemma \ref{Utleq1lemma}, there exists precisely one element in $U_t$ for each
$t \in T^4$, which we will denote by $u_t$.
\item
If $t \in T_s^3$, we put $\omega_t := tsrtst$ and $\pi_t := trs$.
\item
If $t \in T_r^3$, we put $\omega_t := srstrsrt$ and $\pi_t := srsrstr$.
\item
If $t \in T_s^4$, we put $\omega_t := rsrsr\omega_2$ and $\pi_t := rsrsr\omega_2\omega_1tut$, where $u:=u_t$ and $\omega_1$, $\omega_2$ are as in Subsection \ref{omega12subsection}.
\item
If $t \in T_r^4$, we put $\omega_t := r\overline{\omega}_2$ and $\pi_t := r\overline{\omega}_2\overline{\omega}_1utu$, where $u:=u_t$ and $\overline{\omega}_1$, $\overline{\omega}_2$ are as in Subsection \ref{SSrelabel}.
\item
For $t \in T^3$, we put $K_t := J_t$ and
for $t \in T^4$, we put $K_t := J_t  \cup \{u_t\}$.
\item
We put $\widehat{T} := T \cup \{ \infty \}$, $K_{\infty}= J_{\infty} := J$, 
$\omega_{\infty} := 1_W$ and $\pi_{\infty} := srs$.

\item
Finally, for $t\in \widehat{T}$, we put $K_t^{\kdef} = K_t \cup K_t^{\perp}$.
\end{itemize}

\noindent
Let $t \in \widehat{T}$ and if $t\neq \infty$, suppose it is tame. We define
$\delta_t\co K_t^{\kdef} \to \langle K_t^{\kdef} \rangle$ by
$\delta_t(r) = rsr, \delta_t(s) = s$,
$\delta_t(t) = \omega_t t \omega_t^{-1}$ (for $t\neq \infty$),
$\delta_t(u_t) := u_t$ for $t \in T^4$
and $\delta_t(x) := x$ for all $x \in K_t^{\perp}$.

\begin{prop} 
$\delta_t$ is an $(r,s,srs)$-deformation of $K_t^{\kdef}$.                           
\end{prop}

\begin{proof}
This is a consequence of Lemma \ref{H3basiclem}, Corollary \ref{H3basiccor} and Proposition \ref{H4basicdefo} together with its relabeled version.
\end{proof}

\begin{definition}
We call $\delta_t$ the \emph{standard deformation} of $K_t^{\kdef}$.
\end{definition}

\subsection{Tame angle-deformations}

\begin{definition}
Let $K$ be a subset of $S$ containing $J$ and let
$\delta\co K \to \langle K \rangle$ be 
an $(r,s,srs)$-deformation of $K$. Then we call
$\delta$ \emph{tame} if for each $t \in T \cap K$ which is
tame 
in $K$, there exists
an element $w_t  \in \langle K \rangle$ such that
$\delta(x) = w_t \delta_t(x) w_t^{-1}$ for all $x \in K_t^{\kdef}$.
\end{definition}

\noindent
The goal of this section is to prove the following result.

\begin{theorem} \label{mainthm}
There exists a tame $(r,s,srs)$-deformation of $S$. 
\end{theorem}

\subsection{The tame case}

Throughout this subsection, we assume the following.

\begin{itemize}
\item[(TAME)] {\it All elements in $T$ are tame.}
\end{itemize}

\noindent 
For each $t \in \widehat{T}$, let $\delta_t\co K_t^{\kdef} \to \langle K_t^{\kdef}
\rangle$ be the standard deformation.
\\
We put $\widehat{J} := J \cup T \cup J^{\perp}$ and we
define $\hat{\delta}\co \widehat{J} \to \langle \thinspace\widehat{J}\thinspace \rangle$
by $\hat{\delta} \mid_{K_t^{\kdef}} := \delta_t$ for each 
$t \in \widehat{T}$
and $\hat{\delta} \mid_{J^{\perp}} := \id_{J^{\perp}}$.
Note that $\hat{\delta}$ is well-defined.
\\
Our first goal is to prove the following proposition.

\begin{prop} \label{hatprop}
$\hat{\delta}$ is a tame $(r,s,srs)$-deformation of 
$(\langle\thinspace \widehat{J} \thinspace\rangle,\widehat{J}\thinspace)$.
\end{prop}

\begin{lemma} \label{commlemma}
Let $t \in T$ and $y \in J^{\perp}$ such that $o(ty) < \infty$.
Then, $y$ commutes with $\omega_t$ and $\pi_t$.
\end{lemma}

\begin{proof}
If $t\in T^3$ or if $y\neq u_t$, then the tameness of $t$ and (DE1) imply that $y\in J_t^{\perp}$, and hence $y\in K_t^{\perp}$ by Lemma \ref{lemma9.13} and we are done. If $y = u_t$, then the result follows from Lemma \ref{H4basiclem} a).
\end{proof}

\begin{lemma} \label{hatorderlemma}
Let $E := \{ x,y \}$ be an edge of $\widehat{J}$ different
from $J$. Then, there is an element 
$w_E \in \bigcup_{t \in \widehat{T}} \langle J_t \rangle$
such that $\hat{\delta}(E) = E^{w_E}$.
\end{lemma}

\begin{proof}
If $E$ is contained in $J_t$ for some $t \in T$, then
there exists an element $w_E \in \langle J_t \rangle$
such that $\hat{\delta}(E) = E^{w_E}$. This follows from
Lemma \ref{H3basiclem} and Corollary \ref{H3basiccor}.

If $E$ is contained in $J^{\perp} \cup \{ s \}$, then
$\hat{\delta}(E) = E^{1_W}$. Hence, the case $s \in E$ is settled completely. 

Suppose now $x=r$. In this case, we may assume $y \in J^{\perp}$ because
the case $y \in T$ is already covered above. 
For all $y \in J^{\perp}$, we have $ysrs =srsy$ and therefore
$\hat{\delta}(E) = E^{srs}$.

By Lemma \ref{noedgesinT}, it remains to consider the
case where $x \in T$ and $y \in J^{\perp}$. Set $x=t$. As
$\{ x,y \}$ is an edge, it follows from Lemma \ref{commlemma}
that $y$ commutes with $\omega_t$. Hence we have
$\hat{\delta}(E) = E^{\omega_t}$ in this case.
\end{proof}

\noindent
{\bf Proof of Proposition \ref{hatprop}:}
It is readily verified that $\hat{\delta}(\thinspace\widehat{J}\thinspace)$
generates $\langle\thinspace \widehat{J} \thinspace\rangle$ and by
Lemma \ref{hatorderlemma} and Proposition \ref{autoprop} (with ${\bf K} = \{J_t \ | \ t\in \widehat{T}\}$),
it follows that $\hat{\delta}$ extends to an automorphism
$\hat{\alpha}$ of $\langle\thinspace \widehat{J} \thinspace\rangle$, which
implies in particular that $\hat{\delta}(\thinspace\widehat{J}\thinspace)$
is a Coxeter generating set of $\langle\thinspace \widehat{J}\thinspace \rangle$. Using Lemma \ref{hatorderlemma}, it is now straightforward to check that $\hat{\delta}$ satisfies Properties AD1)-AD4).
The tameness of $\hat{\delta}$ is a consequence of its
definition.
This concludes the proof of Proposition \ref{hatprop}. \qedsymbol

\noindent
Let $L$ be a  $J$-component and define the set
$T_L \subseteq T$ as before. Since we assume
that $J$ is a $\Delta$-edge, we know by Corollary \ref{Jtflexcor2} that $|T_L| \leq 1$.
We define $t(L)$ as in Subsection \ref{t(L)subsection}. Moreover, we put $J_L := K_{t(L)}$,
$K_L := J_L \cup J^{\perp} \cup L$ and $M_L := K_L \cup T$.
\\
Let $\Pi(L)$ be the set of $L$-free vertices of $J$;
since $J$ is flexible (by Lemma \ref{Jflexlemma}), we know that $\Pi(L) \neq \emptyset$.

For each $J$-component $L$, we define $\gamma_L \in \langle J_L \rangle$
as follows.
\\
If $t(L) = \infty$ and $r \in \Pi(L)$, we put $\gamma_L := 1_W$.
\\
If $t(L) = \infty$ and $ \Pi(L) = \{ s \}$, we put $\gamma_L := srs$.
\\
If $t(L) \in T_s$ and $r \in \Pi(L)$, we put $a_L := r$ and $\gamma_L := \omega_t$.
\\
If $t(L) \in T_s$ and $ \Pi(L) = \{ s \}$, we put $a_L:= s$ and $\gamma_L := \pi_t$.
\\
If $t(L) \in T_r$ and $s \in \Pi(L)$, we put $a_L:= s$ and $\gamma_L := \omega_t$.
\\
If $t(L) \in T_r$ and  $ \Pi(L) = \{ r \}$, we put $a_L := r$ and $\gamma_L := \pi_t$.
\\
Finally, we define $\delta_L\co K_L \to \langle K_L \rangle$
by $\delta_L \mid_{J_L} := \delta_{t(L)} \mid_{J_L}$, 
$\delta_L \mid_{J^{\perp}} =
\id_{J^{\perp}}$ and $\delta_L(x) := \gamma_L x \gamma_L^{-1}$
for all $x \in L$. Note that $\delta_L$ is well-defined.

\begin{lemma}
Let $L$ be a $J$-component with $t:= t(L) \neq \infty$.
Then $K_t$ is an $a_L$-special subset of $K_L$.
\end{lemma}
                                                                                                                             
\begin{proof} This is a consequence of Lemma \ref{lemma9.14}. 
\end{proof}
     
\begin{lemma} \label{Lcomplemma}
Let $L$ be a $J$-component. Then $\delta_L$ is
an $(r,s,srs)$-deformation of $K_L$.
\end{lemma}

\begin{proof} This is a consequence of the previous lemma and Propositions \ref{defexiprop1}, \ref{sTWprop} and \ref{relabelprop} applied to the Coxeter system $(\langle K_L \rangle, K_L)$ if $t(L)\neq\infty$, and of Proposition \ref{rk2cpdefprop} applied to the same Coxeter system otherwise.
\end{proof}

\begin{prop} \label{Lcompprop} 
Let $L$ be a $J$-component. Define
$\hat{\delta}_L\co M_L \to \langle M_L \rangle$
by $\hat{\delta}_L \mid_{K_L} := \delta_L$ and
$\hat{\delta}_L \mid_{\widehat{J}} := \hat{\delta}$.
Then $\hat{\delta}_L$ is an $(r,s,srs)$-deformation of $M_L$.
\end{prop}

\begin{proof} Note first that $\hat{\delta}_L$ is well-defined. 
By the previous lemma, $\delta_L$ is 
an $(r,s,srs)$-deformation of $K_L$ and by 
Proposition \ref{hatprop}, $\hat{\delta}$ is
an $(r,s,srs)$-deformation of $\widehat{J}$. As 
$K_L \setminus \widehat{J} = L$ and 
$\widehat{J} \setminus K_L = T \setminus T_L$, all edges
of $M_L$ are contained in at least one of the two sets.
Now, as $\delta$ restricted to $M_L\cap K_L = J_L \cup J^{\perp}$ is an $(r,s,rsr)$-deformation of $M_L\cap K_L$, Proposition \ref{amalprop} finishes the proof.
\end{proof}

\begin{theorem} \label{Tamethm}
Let $\delta\co S \to W$ be the mapping defined
by $\delta \mid_{M_L} := \hat{\delta}_L$ for each $J$-component
$L$. Then $\delta$ is a tame
$(r,s,srs)$-deformation of $(W,S)$. 
\end{theorem}
 
\begin{proof} 
Note that for two different $J$-components $L$ and $L'$, we have $M_L\cap M_{L'} = T\cup J^{\perp}\cup J$, which is independent of $L$ and $L'$. Moreover, $\delta$ restricted to $T\cup J^{\perp}\cup J$ is an $(r,s,rsr)$-deformation of $T\cup J^{\perp}\cup J$.
The claim now follows by induction on the
number of $J$-components using Propositions \ref{amalprop}
and \ref{Lcompprop}, the tameness being a consequence of Proposition \ref{hatprop}.
\end{proof}
                                                                                                                            
\subsection{Proof of Theorem \ref{mainthm}} 

The theorem will be proved by induction on $\deg(S)$.
If $\deg(S) = 0$, all elements in $T$ are tame and
we are done by Theorem \ref{Tamethm}. 
Suppose now that the degree of $S$ is at least 1.
Then there exists a wild $t \in T$, which we fix throughout
this subsection.

\noindent 
For each $u \in \widehat{U}_t$, we define
the sets $V_u,W_u,X_u,Y_u$ and $Z_u$ as in Subsection
\ref{WXYZsubsection}.
For $u \in U_t$, we put $\tau_u := \tau$ where $\tau$
is defined as in Subsection \ref{omega12subsection} and
$\tau_{\infty} := 1_W$.

Let $u \in \widehat{U}_t$.
By Lemma \ref{deglemma}, we know that $\deg(Y_u) < \deg(S)$.
Thus, we know by induction
that there is  a tame $(r,s,srs)$-deformation $\theta_u$
of $Y_u$. Again by Lemma \ref{deglemma},
$t$ is tame in $Y_u$ and if we define $K_t^{\kdef}$ as
in Subsection \ref{sstameconv} with respect to $Y_u$, we
have $W_u = K_t^{\kdef}$. 
Hence, there is an $w_u \in \langle Y_u \rangle$
such that $\Int(w_u) \circ \theta_u \mid_{W_u}$ is the
standard deformation of $W_u$. We put $\Theta_u:= \Int(w_u) \circ \theta_u$.
The discussion above yields the following.

\begin{lemma}
For each $u \in \widehat{U}_t$, there exists a tame
$(r,s,srs)$ deformation $\Theta_u$ of $Y_u$ such that
$\Theta_u \mid_{W_u}$ is the standard deformation of
$W_u$.
\end{lemma}

\noindent
For each $u \in \widehat{U}_t$, let $\Theta_u$ be as in the
previous lemma and 
put $\delta_u:= \Int(\tau_u) \circ \Theta_u$. 

\begin{lemma} \label{lemma1042} 
For each $u \in \widehat{U}_t$ the mapping $\delta_u\co Y_u \to 
\langle Y_u \rangle$ is a tame $(r,s,srs)$-deformation of $Y_u$. 
Moreover, we have $\delta_u \mid_{Y_u \cap Y_{\infty}}
= \delta_{\infty} \mid_{Y_u \cap Y_{\infty}}$. In particular,
there exists an $(r,s,srs)$-deformation $\hat{\delta}_u$ of $Z_u$
such that $\hat{\delta}_u \mid_{Y_u} = \delta_u$ and
$\hat{\delta}_u \mid_{Y_{\infty}} = \delta_{\infty}$.
\end{lemma}

\begin{proof}
The first assertion of the lemma is clear, because $\Theta_u$
is tame and $\tau_u \in \langle Y_u \rangle$. 

The second assertion is trivial for $u=\infty$, so we may assume $u\in U_t$. First remark that $Y_u \cap Y_{\infty} = J_t \cup J_{t,u}^{\perp}$ by Lemma \ref{lemma923}.
Since $\Theta_u \mid_{W_u}$ is the standard deformation
and as $\tau_u \in \langle J_{t,u} \rangle$ 
commutes with all elements
in $J_{t,u}^{\perp}$ and with $rsr$ and $s$ (by Lemma \ref{H4basiclem} c)),
it follows that 
$\delta_u \mid_{J \cup J_{t,u}^{\perp}} = \delta_{\infty} \mid_{J \cup J_{t,u}^{\perp}}$. Thus, it remains only to check whether
$\delta_u(t) = \delta_{\infty}(t)$; but this is also a consequence
of Lemma \ref{H4basiclem} c). This concludes the proof of the second
assertion.

The last assertion is a consequence of the second, Lemma \ref{lemma923}
 and Proposition \ref{amalprop}.
\end{proof}

\begin{lemma}
There exists an $(r,s,srs)$-deformation $\delta$ of
$S$ such that $\delta \mid_{Y_u} = \delta_u$ for each
$u \in \widehat{U}_t$.
\end{lemma}
 
\begin{proof} As $t$ is assumed to be wild, we have
$|U_t| \geq 1$. We prove the lemma by induction
on $|U_t|$. If $|U_t|=1$ and if $u$ denotes the unique
element in $U_t$, then $S=Z_u$ and we are done 
by the previous lemma.

Suppose now $|U_t| > 1$ and let $u \in U_t$.
Put $C_u := \bigcup_{u \neq u' \in U_t} Z_{u'}$.
Note first that $C_u \cap Z_u = Y_{\infty}$ and that
each edge in $S$ is contained in $C_u$ or in $Z_u$.
By induction, there exists an $(r,s,srs)$-deformation $\delta_u'$
of $C_u$  such that $\delta_u' \mid_{Y_{u'}} = \delta_{u'}$
for each $u' \in \widehat{U}_t$ different from $u$.
By the previous lemma, there exists an $(r,s,srs)$-deformation
$\hat{\delta}_u$
of $Z_u$ such that $\hat{\delta}_u \mid_{Y_a} = \delta_a$
for $a \in \{ u, \infty \}$. 
Now Proposition \ref{amalprop} yields the existence
of $\delta$.  
\end{proof}

\noindent
{\bf Conclusion of the Proof of Theorem \ref{mainthm}:}
The previous lemma yields the existence of an $(r,s,srs)$-deformation
$\delta$ of $S$ such that $\delta \mid_{Y_u} = \delta_u$
for each $u  \in \widehat{U}_t$. It remains to show that
$\delta$ is tame. Let $t' \in T$ be tame in $S$. Since $t$ is assumed
to be wild, we have $t' \neq t$. By Lemmas \ref{lemma922} and 
\ref{lemma9.13}, there is an $u \in \widehat{U}_t$ such that
$K_{t'}^{\kdef}$ is contained in $Y_u$. By Lemma \ref{lemma1042},
we know that $\delta_u$ is a tame $(r,s,srs)$-deformation.
Hence there exists an element $v \in \langle Y_u \rangle$
such that $\Int(v) \circ \delta_u \mid_{K_{t'}^{\kdef}}$
is the standard deformation of $K_{t'}^{\kdef}$. As
$\delta \mid_{Y_u} = \delta_u$, it follows that 
$\Int(v) \circ \delta \mid_{K_{t'}^{\kdef}}$ is the
standard deformation of $K_{t'}^{\kdef}$. Hence
$\delta$ is tame.  \qedsymbol

\section{Proof of Theorem \ref{thm2}}
\label{sec11}
Let $(W,R)$ be a Coxeter system and let $S\subseteq R^W$ be a Coxeter generating set which is not sharp-angled. Suppose $S$ contains $k\geq 1$ edges which are not sharp-angled and choose one of them. Call it $J$. By Theorem \ref{thm1}, we can assume that $J=\{r,s\}$ with $o(rs)=5$. By Proposition \ref{deftheorem}, $J$ is a $\Delta$-edge. Hence, by Theorem \ref{mainthm}, there exists a $J$-deformation $\delta$ of $S$ sending $J$ onto $\{rsr,s\}$. Hence, by Lemma \ref{rk2deflem}, $\delta(J)$ is a sharp-angled edge of $\delta(S)$. Let now $J'$ be an edge of $S$ different from $J$. Then $\delta(J')$ is $W$-conjugate
to $J'$ by Property AD4) of $\delta$; in particular, $\delta(J')$ is sharp-angled
if and only if $J'$ is sharp-angled. Hence the number of edges in $\delta(S)$
which are not sharp-angled is $k-1$. Thus the statement follows by an
obvious induction on the number of edges of $S$
which are not sharp-angled. \qedsymbol

\bibliographystyle{amssort}

\end{document}